\documentclass[EJP]{ejpecp}  

\usepackage{amsmath}
\usepackage{amssymb}
\usepackage{amsthm}
\usepackage{psfrag}
\usepackage{verbatim}
\usepackage{color}

\SHORTTITLE{Dynamical freezing} 

\TITLE{Dynamical Freezing in a Spin Glass System with Logarithmic Correlations} 

\AUTHORS{%
  Aser~Cortines\footnote{University of Z\"urich, Switzerland.
    \EMAIL{aser.cortinespeixoto@math.uzh.ch}}
  \and
  Julian~Gold\footnote{Northwestern University, United States of America.
    \EMAIL{gold@math.northwestern.edu}}
  \and
  Oren~Louidor\footnote{Technion, Israel.
    \EMAIL{oren.louidor@gmail.com}}
    }


\AMSSUBJ{60K37; 82C44; 60G57; 60G70; 60G15; 82D30} 

\SUBMITTED{11/26/2017} 
\ACCEPTED{5/23/2018} 

\ARXIVID{1711.01360} 

\VOLUME{23}
\YEAR{2018}
\PAPERNUM{59}
\DOI{10.1214/YY-TN}

\ABSTRACT{We consider a continuous time random walk on the two-dimensional discrete torus, whose motion is governed by the discrete Gaussian free field on the corresponding box acting as a potential. More precisely, at any vertex the walk waits an exponentially distributed time with mean given by the exponential of the field and then jumps to one of its neighbors, chosen uniformly at random. We prove that throughout the low-temperature regime and at in-equilibrium timescales, the process admits a scaling limit as a spatial K-process driven by a random trapping landscape, which is explicitly related to the limiting extremal process of the field. Alternatively, the limiting process is a supercritical Liouville Brownian motion with respect to the continuum Gaussian free field on the box. This demonstrates rigorously and for the first time, as far as we know, a dynamical freezing in a spin glass system with logarithmically correlated energy levels. }

\newtheorem{lem}[theorem]{Lemma}
\newtheorem{prop}[theorem]{Proposition}

\newtheorem{mainthm}{Theorem}

\newcommand{\bbE}{\mathbb{E}}

\newcommand{\bbP}{\mathbb{P}}
\newcommand{\bbN}{\mathbb{N}}
\newcommand{\bbZ}{\mathbb{Z}}
\newcommand{\bbR}{\mathbb{R}}

\newcommand{\eps}{\epsilon}

\newcommand{\vart}{\vartheta}
\newcommand{\ol}{\overline}

\newcommand{\wh}{\widehat}

\newcommand{\cA}{\mathcal A}

\newcommand{\cC}{\mathcal C}

\newcommand{\cS}{\mathcal S}

\newcommand{\cF}{\mathcal F}
\newcommand{\cM}{\mathcal M}

\newcommand{\rmc}{{\rm c}}

\newcommand{\rmd}{{\rm d}}
\newcommand{\rme}{{\rm e}}
\newcommand{\rmB}{{\rm B}}

\newcommand{\rmL}{{\rm L}}

\newcommand{\wt}{\widetilde}

\newcommand{\sfS}{{\text{\rm\textsf{S}}}}
\newcommand{\sfX}{{\text{\rm\textsf{X}}}}
\newcommand{\sfR}{{\text{\rm\textsf{R}}}}
\newcommand{\sfJ}{{\text{\rm\textsf{J}}}}
\newcommand{\sfE}{{\text{\rm\textsf{E}}}}
\newcommand{\sfH}{{\text{\rm\textsf{H}}}}
\newcommand{\sfP}{{\text{\rm\textsf{P}}}}
\newcommand{\sfG}{{\text{\rm\textsf{G}}}}
\newcommand{\sfL}{{\text{\rm\textsf{L}}}}
\newcommand{\sfT}{{\text{\rm\textsf{T}}}}

\newcommand{\e}{\epsilon}
\newcommand{\pa}{\partial}

\newcommand{\1}{\textbf{1}}
\newcommand{\mft}{\mathfrak{t}}
\newcommand{\Leb}{{\text{\rm\text{Leb}}}}

\newcommand{\note}[1]{}

\begin{document}



\section{Introduction}

\subsection{Setup and main result}
\label{ss:SetupAndResults}
Let $V_N := [0, N)^2 \cap \bbZ^2$ be the discrete box of side length $N$ and consider the random field $h_N \equiv (h_{N,x})_{x \in V_N}$ having the law of the discrete Gaussian free field (DGFF) on $V_N$ with zero boundary conditions. That is, $h_N$ is a centered Gaussian with covariance given by $\sfG_{V_N}$, the discrete Green function of a simple random walk on $\bbZ^2$ killed upon exiting $V_N$. 

On the same probability space, define a process $X_N \equiv (X_N(t) : t \geq 0)$ taking values in the two-dimensional discrete torus $V_N^* = \bbZ^2/N\bbZ^2$, whose vertices we identify with those of $V_N$. The law of $X_N$ is given {\em conditionally on $h_N$} to be Markovian with transition rates:
\begin{equation}
q^{(\beta)}_N(x,y) := \left \{ \begin{array}{ll}
	\tfrac{1}{4}\rme^{-\beta h_{N,x}}	& \text{if } d_N^*(x,y)=1 \,, \\
	0					& \text{otherwise}\,,
	\end{array} \right.
\label{eq:1.1_dynamics}
\end{equation}
where $d_N^*$ is the metric on $V^*_N$, which is naturally induced by the Euclidean norm $\|\cdot\|$ on $\bbZ^2$.
Evidently, $X_N$ is a continuous-time symmetric random walk on $V_N^*$ with mean holding time $\rme^{\beta h_{N,x}}$ at vertex $x$. Moreover, it is reversible with respect to its unique stationary distribution $\wh{\mu}^{(\beta)}_N$ given by
\begin{equation}
\label{eq:1.9}
\wh{\mu}^{(\beta)}_N := \frac{\mu^{(\beta)}_N}{\mu^{(\beta)}_N(V_N)}
\quad, \qquad
\mu^{(\beta)}_N := \sum_{x \in V_N} \rme^{\beta h_{N,x}} \delta_{x} \,.
\end{equation}

The aim of the present work is to derive a large $N$ distributional limit for the process $X_N$ when time and space are scaled appropriately, and we therefore start by constructing the limiting object. The first ingredient in this construction is a spatial version of the {\em $K$-process}, which was introduced by 
Fontes and Mathieu in~\cite{fontes2008k} (see Subsection~\ref{sss:KProcess} for further historical comments). 

To define this process, let $\tau \equiv (\tau_k)_{k=1}^\infty$ be a decreasing sequence of positive values, thought of as {\em trap depths} (the ``trapping'' terminology will become apparent in the next subsection). For $k \geq 1$, let $A_k \equiv (A_k(u) : {u \geq 0})$ be independent Poisson processes with rate $1$,
and let $\sigma_j^{(k)}$ denote the time of the $j$-th jump of $A_k$, with the convention that $\sigma_{0}^{(k)} := 0$. For $u \geq 0$, define the \emph{clock process} $T^{(\tau)} \equiv (T^{(\tau)}(u) : u \geq 0)$ associated with $\tau$ by,
\begin{align}
T(u) \equiv T^{(\tau)}(u) := \sum_{k=1}^\infty \tau_k \sum_{j=1}^{A_k(u)} e_j^{(k)} \,,
\label{eq:clock}
\end{align}
where $(e_j^{(k)})_{j,k = 1}^\infty$ are i.i.d. exponentially distributed random variables with rate $1$, which are independent of the $A_k$'s. 

By taking expectation, one finds that $T^{(\tau)}(u) < \infty$ for all $u \geq 0$ almost surely, provided that $\sum_k \tau_k < \infty$. In this case, we may construct a process $K^{(\tau)} \equiv (K^{(\tau)}(t) : t \geq 0)$ taking values in $\bbN^* := \bbN \cup \{\infty\}$ via
\begin{equation}
K(t) \equiv K^{(\tau)}(t) := \left\{
	\begin{array}{ll}
	k			& \text{if } t \in \big[T(\sigma_{j}^{(k)}-)\,,\,\, T(\sigma_{j}^{(k)})\big) 
					\text{ for some } k \geq 1 \text{ and } j \geq 1 \,, \\
	\infty		& \text{otherwise.}
	\end{array} \right.
\label{eq:K_def}
\end{equation}
We refer to $K^{(\tau)}$ as the K-process associated with depths $\tau$. 

Given an additional sequence $\xi \equiv (\xi_k)_{k=1}^\infty$ of points in $V := [0,1]^2$, we may also define the process $Y^{(\xi, \tau)} \equiv (Y^{(\xi, \tau)}(t) : t \geq 0)$ taking values in $V^* := V \cup \{\infty\}$ via
\begin{equation}
\label{eq:spatial_K}
 Y(t) \equiv Y^{(\xi,\tau)}(t) := 
 \begin{cases}
 	\xi_{K(t)} & \text{if } K(t) \neq \infty\,, \\
 	\infty & \text{if } K(t) = \infty \,.
 \end{cases}
\end{equation}
We refer to $\xi$ as a sequence of {\em trap locations} and to $Y^{(\xi, \tau)}$ as 
the {\em spatial K-process} associated with {\em trapping landscape} $(\xi, \tau)$.

Next, we recall some results from the extreme value theory for the DGFF. To this end, define the {\em structured extremal process} of $h_N$ as the point process on $V \times \bbR \times [0,\infty)^{\bbZ^2}$ given by
\begin{equation}
\label{e:1.5}
\eta_{N,r} := \sum_{x \in V_N} \delta_{x/N} \otimes \delta_{h_{N,x} - m_N} \otimes
	\delta_{(h_{N,x} - h_{N,x+y})_{y \in \rmB_r}} \1_{\left\{h_{N,x} \geq h_{N,x+y} \,:\, y \in \rmB_r \right\}} \,,
\end{equation}
where $m_N := 2\sqrt{g}\log N - \tfrac34 \sqrt{g} \log \log N$ is an appropriate centering sequence and $\rmB_r = \rmB_r(0)$ is the open Euclidean ball of radius $r$ around $0$ in $\bbZ^2$. Setting $g := 2/\pi$ and $\alpha = \sqrt{2\pi}$, it was shown in~\cite[Theorem 2.1]{biskup2016full} that 
there is a probability measure $\nu$ on $[0,\infty)^{\bbZ^2}$ and a random measure $Z$ on $V$, satisfying $|Z| := Z(V) \in (0,\infty)$ almost surely, such that with $\wh{Z} = Z/|Z|$, we have
\begin{equation}
\label{eq:thm1}
\eta_{N, r} \underset{\substack{N \to \infty\\ r \to \infty}} \Longrightarrow \eta \sim \text{\rm PPP}\bigl(\wh{Z}(\rmd x)\otimes |Z| \rme^{-\alpha h} \rmd h \otimes \nu(\rmd \phi)\bigr) \,.
\end{equation}
Above PPP stands for Poisson point process, and the law of $\eta$  
should be interpreted as given conditionally on $Z$, which is assumed to be defined on the same probability space as $\eta$ itself. Also, here and after, whenever we write $F(p,q) \to F(p_0, q_0)$ in the limit when $p \to p_0$ followed by $q \to q_0$, we mean that $\lim_{q \to q_0} \limsup_{p \to p_0} d\big(F(p,q), F(p_0, q_0)\big) = 0$, where $d$ is an appropriate metric on the range of $F$. 

The second ingredient in the construction of the limiting object is a point process which we denote, for $\beta > \alpha$, by $\chi \equiv \chi^{(\beta)}$. It is defined on the same probability space as the measure $Z$ and its law is given conditionally on $Z$ as
\begin{equation}
\label{eq:chi}
\chi^{(\beta)} \sim \text{PPP} \big(\wh{Z}(\rmd z) \otimes \kappa_{\beta} |Z| t^{-1-\alpha/\beta}  \rmd t \big) \,,
\end{equation} 
where
\begin{equation}
\label{eq:kappa}
\kappa_\beta
:= \frac{1}{\beta} \int \Big(\sum_{y \in \bbZ^2} \rme^{-\beta \omega_y} \Big)^{\alpha/\beta} \nu(\rmd  \omega) 
\,.
\end{equation}
The integral above, which involves the distribution $\nu$ from the law of $\eta$, was shown to be finite in~\cite[Theorem 2.6]{biskup2016full}. Together with the stated properties of $Z$, this ensures that $\chi^{(\beta)}$ is a well-defined random Radon point measure on $V \times (0, \infty)$. Furthermore, the almost-sure integrability of $t \mapsto (t \wedge 1)$ under the second component of the conditional intensity measure assures that we can enumerate the atoms of $\chi^{(\beta)}$ in descending order of their second coordinate, and that if $(\xi^{(\beta)}, \tau^{(\beta)}) \equiv \big((\xi^{(\beta)}_k, \tau^{(\beta)}_k)\big)_{k=1}^\infty$ denotes such an enumeration, then $\sum_k \tau^{(\beta)}_k < \infty$ almost surely (see Lemma~\ref{lem:finite_tau}). We refer to $(\xi^{(\beta)}, \tau^{(\beta)})$ as the {\em limiting trapping landscape of the DGFF}. 

Combining both ingredients, we can now construct the process $Y^{(\beta)} \equiv \big(Y^{(\beta)}(t) : t \geq 0\big)$ on the same probability space as that of $\chi^{(\beta)}$ by specifying its conditional law as
\begin{equation}
Y^{(\beta)} \,|\, \chi^{(\beta)} \overset{d}{=} Y^{(\xi^{(\beta)}, \tau^{(\beta)})} \,.
\end{equation}
The process $Y^{(\beta)}$, henceforth called the {\em $\chi$-driven spatial K-process}, is the desired limiting object. 

For what follows, we endow $V^*$ with the metric $d^*$ which agrees with the Euclidean distance on $V$ and puts $\infty$ at a distance $1$ from all other points (the choice of $1$ is arbitrary and any positive number would yield an equivalent metric). For $\mft > 0$, the space $\rmL([0,\mft], V^*)$ consists of all measurable functions from $[0,\mft]$ to $V^*$, and is equipped with the metric $\|f-g\|_{\rmL([0,\mft], V^*)} := \int_0^{\mft} d^*(f(s),g(s)) \rmd s$ for 
$f,g \in \rmL([0,\mft], V^*)$. Observe that under $\|\cdot\|_{\rmL([0,\mft], V^*)}$ (which is not a norm, despite the notation) this space is complete and separable, and that this metric generates the topology of convergence in measure on functions from $[0, \mft]$ to $V^*$. Here, the interval $[0,\mft]$ is implicitly equipped with Lebesgue measure. The next theorem is the principle result of this work.

\begin{mainthm}
\label{thm:A}
Let $\beta > \alpha$. Then with $s_N := g N^{2\sqrt{g} \beta}(\log N)^{1-3\sqrt{g}\beta/4}$ and for any $\mathfrak{t} > 0$,
\begin{equation}
\label{e:1.8}
\Big(\tfrac{1}{N} X_N(s_Nt) : t \in [0,\mathfrak{t}] \Big) \Longrightarrow \Big(Y^{(\beta)}(t) : t \in [0,\mathfrak{t}]\Big) 
\quad \text{as } \ N \to \infty \,,
\end{equation}
where the above weak convergence is that of random functions in $\rmL\big([0,\mathfrak{t}], V^*\big)$.
\end{mainthm}

An equivalent formulation of this theorem, which is just an explicit rewrite of~\eqref{e:1.8}, is as follows: for any $\epsilon > 0$ if $N$ is large enough, then there is a coupling between $X_N$ and $Y^{(\beta)}$ so that with probability at least $1-\epsilon$, the set of times $t \in [0,\mft]$ with $d^*\big(\tfrac{1}{N} X_N (s_N t),\, Y^{(\beta)}(t)\big) > \epsilon$ has Lebesgue measure at most $\epsilon$.

\subsection{Interpreting the limiting process}
Theorem~\ref{thm:A} states that on timescales of the order $s_N$ and when space is scaled by $N$, the process $X_N$ tends to be close in law to the process $Y^{(\beta)}$. Understanding the large $N$ behavior of $X_N$ therefore requires an understanding of the statistics of $Y^{(\beta)}$ and in turn the probabilistic features of its two ingredients: the spatial K-process $Y^{(\xi, \tau)}$ and the limiting trapping landscape of the DGFF $(\xi^{(\beta)}, \tau^{(\beta)})$.

\subsubsection{The spatial K-process}
\label{sss:KProcess}
The process $Y^{(\xi, \tau)}$ is a spatial version of the $\bbN \cup \{\infty\}$-valued K-process of Fontes and Mathieu, which was introduced in~\cite{fontes2008k} to describe scaling limits of dynamics in effective trap models. It is named after Kolmogorov~\cite{kolmogorov}, who considered a process similar to $K^{(\tau)}$, albeit with the clock process $T(u)$ in~\eqref{eq:clock} modified by a linear term. Endowing $\bbN \cup \{\infty\}$ with any metric which makes it compact, Fontes and Mathieu show that $K^{(\tau)}$ is strongly Markovian on $\bbN \cup \{\infty\}$ and can be defined to have c\`{a}dl\`{a}g sample paths.

The evolution of $K^{(\tau)}$ is as follows. When at state $k \in \bbN$, the process waits an exponentially distributed time with mean $\tau_k$ and then jumps to $\infty$. At $\infty$ it spends $0$ time (that is the exit time is almost surely $0$), after which it jumps to a new state $k' \in \bbN$ (possibly $k$ again). Following a visit to $\infty$, the hitting time of any finite subset of states in $\bbN$ is finite almost surely and the entrance distribution to any such subset is uniform. Finally, the set of times at which the process is at $\infty$ has Lebesgue measure $0$. In the terminology of Markov processes, $\infty$ is said to be an {\em unstable} and {\em fictitious} state.

Thinking of $\bbN$ as indexing countably many traps, with trap $k$ having depth $\tau_k$, the process $K^{(\tau)}$ represents a dynamics of ``uniform'' trap hopping. That is, it evolves by hopping from one trap to another, getting stuck at a trap for a time with mean corresponding to its depth and then jumping to the next trap, which is chosen (formally) uniformly at random. Facilitating such a jump is the state $\infty$, in which the process spends an infinitesimal amount of time, and which essentially captures this transitional period. 
 
In analog, if traps are thought of as placed on $V$, such that trap $k$ is at location $\xi_k$, then the spatial K-process with trapping landscape $(\xi, \tau)$, as given in~\eqref{eq:spatial_K}, describes the current location (instead of trap index) of a dynamics which is defined exactly as before. Observe that the process is still Markovian and that $\infty$, which still represents a state of transition from one trap to another, remains unstable and fictitious.

\subsubsection{The limiting trapping landscape of the DGFF}
\label{sss:DGFFLandscape}

In view of the law of $\chi^{(\beta)}$, the depths of the traps in the trapping landscape are distributed as the atoms, in descending order, of a Poisson point process on $(0,\infty)$ with the Fr\'{e}chet intensity measure $t^{-\alpha/\beta - 1} \rmd t$, all multiplied by the global random factor $(\kappa_{\beta} |Z|)^{\beta/\alpha}$. 
The corresponding trap locations are then i.i.d., and drawn independently of the depths from distribution $\wh{Z}$.  Both $|Z|$ and $\wh{Z}$ are determined from the random measure $Z$, which is drawn beforehand. 

To explain the appearance of $\chi^{(\beta)}$ in the definition of the trapping landscape governing $Y^{(\beta)}$, we appeal to Theorem~2.6 from~\cite{biskup2016full}, which states that when $\beta > \alpha$, the measures $\mu_N^{(\beta)}$ and $\wh{\mu}_N^{(\beta)}$ from~\eqref{eq:1.9} admit the following large $N$ distributional scaling limits:
\begin{equation}
\label{e:1.11}
\rme^{-\beta m_N} \mu_N^{(\beta)}(N\cdot) \  \Longrightarrow \ 
\mu^{(\beta)} \! := \textstyle{\sum_{(x,t) \in \chi^{(\beta)}}} \, t \delta_x
\quad, \qquad
\wh{\mu}_N^{(\beta)}(N \cdot) \ \Longrightarrow\  \wh{\mu}^{(\beta)} \! := \! \frac{\mu^{(\beta)}}{\mu^{(\beta)}(V)} \,,
\end{equation}
where both limiting measures are purely discrete.

Thus, the process $\chi^{(\beta)}$ is encoding the atoms (location and mass) of the limiting measure $\mu^{(\beta)}$ and, after normalization by the sum of the masses, also the atoms of the limiting stationary distribution $\wh{\mu}^{(\beta)}$. It follows that in equilibrium $N^{-1} X_N$ is typically found close to one of the locations $\xi^{(\beta)}_k$ with probability proportional to $\tau^{(\beta)}_k$, and it is therefore no surprise that when timescales are tuned properly, this static behavior translates into a hopping dynamics with respect to $\xi^{(\beta)}$ and $\tau^{(\beta)}$.

We remark that when $\beta \leq \alpha$, the measures $\mu_N{(\beta)}$ and $\wh{\mu}^{(\beta)}_N$ still admit scaling limits (with a different normalization), but the limiting measures are no longer discrete (see, for instance, \cite{biskup2016intermediate, rhodes2013gaussian}). In the terminology of spin glasses one says that $\beta = \alpha$ marks the {\em glassy transition} point of the system (see Subsection~\ref{sss:SpinGlasses}).

\subsubsection{The $\chi$-driven spatial K-process}
Altogether, the large $N$ behavior of $N^{-1} X_N (s_N \cdot)$ is that of uniform trap hopping dynamics, with an underlying trapping landscape reflecting (up to a global multiplicative factor) the atoms of the equilibrium distribution of $X_N$. Combining the descriptions above, we see that $X_N$ tends to get stuck in clusters of {\em meta-stable} states where its stationary distribution is exceptionally large. These clusters have diameter $o(N)$  (in fact $O(1)$, see Subsection~\ref{ss:ExtremalPic}), and $X_N$ jumps over or \emph{tunnels} through the remaining vertices in negligible time (see~\cite{beltran2010tunneling} for a  mathematical framework for meta-stability and tunneling, which is aligned with our use of these terms).

\subsection{Motivation and related work}
Let us now discuss the motivation for studying the model in this work. A key feature of the underlying field $h_N$ is its logarithmic correlations. Indeed, known asymptotics for the discrete Green function (see Lemma~\ref{lem:GAsymp}) yield
\begin{equation}
\label{eq:log_cor}
\bbE (h_{N,x} - h_{N,y})^2 = \sfG_N(x,x) + \sfG_N(y,y) - 2 \sfG_N(x,y) =  2g \log \|y-x\| + O(1) \,,
\end{equation}
for $x,y \in V_N$ away from the boundary of $V_N$. Such fields have attracted considerable attention in the past few years, but mainly for their ``static'' (structural) features. Nevertheless, logarithmic correlations are of considerable interest in the dynamical context as well, and we proceed to discuss several possible interpretations of the model which lead to different motivations for considering it.

\subsubsection{Low temperature spin glass dynamics and effective trap models}
\label{sss:SpinGlasses}
One may view $\wh{\mu}^{(\beta)}_N$ from~\eqref{e:1.11} as a spin glass-type Gibbs distribution for a thermodynamical system at temperature $\beta^{-1}$ with energy states $(-h_{N,x} :\: x \in V_N)$. In this context, the process $X_N$ models the Glauber dynamics by which such a system relaxes to its equilibrium state. A conjectured universal feature of many spin glass systems, at least of the mean-field type, is the occurrence of a {\em glassy phase} at low temperature. This phase is characterized by a trapping behavior (also meta-stability or freezing) in the dynamics which reflects a stationary distribution with few isolated and dominant energy states (or clusters of these states).

So far this low-temperature dynamical picture has only been verified mathematically in few instances, with {\em Derrida's random energy model} (REM) and {\em Bouchaud's trap model} (BTM) being the principal ones. The REM is defined similarly to our model, only with $V_N^*$ usually replaced by the hypercube $H_N := \{0,1\}^n$, where $N=2^n$ and, more importantly, where the Gaussian potential $h_N = (h_{N,x} :\: x \in H_N)$ is uncorrelated with variance diverging logarithmically in $N$. In the BTM, which was introduced as an effective trapping model, one considers any underlying graph, but replaces $h_{N,x}$ with $\log \tau_x$, where the $\tau_x$'s are i.i.d. and in the domain of attraction of a stable law with index $1$. 

At short (pre-equilibrium) timescales, both the REM and the BTM on various finite graphs and with $\beta > 1$, have been shown to exhibit {\em aging} (see, for example,~\cite{BBGglauber_1, BBGglauber_2,arous2008arcsine}). The latter is a trapping phenomenon: the evolving dynamics encounters states with progressively lower energy levels, around which it gets stuck for increasingly longer periods of time and thus gradually slows down (or ages). Other cases in which pre-equilibrium timescales were considered include Glauber dynamics for the SK $p$-spin model~\cite{bovier2013convergence} and the BTM on $\bbZ^d$~\cite{BC_scaling_Zd}.

Longer (in-equilibrium) timescales were studied in the BTM with $\beta > 1$ on various finite graphs~\cite{fontes2008k, JLT} and under fairly general graph-topological conditions~\cite{jara2014universality}. In all cases, the dynamics were shown to converge, under proper scaling of time and space, to the trapping dynamics given by the (non-spatial) K-process. As in our case, trap depths are given by the atoms of a Poisson point process with the Fr\'{e}chet intensity, albeit without the random multiplicative factor as in~\eqref{eq:chi}. In-equilibrium timescales were also studied in the case of Glauber dynamics for the generalized REM, see~\cite{fontes2014grem,fontes2014k}.

In view of~\eqref{eq:log_cor} and the discussion in Subsection~\ref{sss:DGFFLandscape}, our work can be seen as demonstrating dynamical freezing for a spin glass system with a logarithmically correlated potential, throughout its glassy phase $(\beta > \alpha)$, and observed at in-equilibrium timescales.  We stress that logarithmic correlations are natural to consider for such systems, as they reflect the conjectured ultra-metric correlation structure of energy states in more realistic spin glass models at low temperature. 

Although we recover the K-process in the limit, as in the case of the BTM (thereby strengthening its position as a universal object for in-equilibrium dynamics of spin glasses in their glassy phase), there are three notable differences compared to the i.i.d. case. First, traps are not single vertices but rather clusters having, essentially, finite diameter (although this diameter disappears in the scaling limit). Second, traps are not scattered uniformly on the underlying domain, but rather according to the non-trivial distribution $\wh{Z}$, which is itself random. Third, there is an overall random multiplicative factor $|Z|$ governing the depths of all traps, resulting in a global random slow-down or speed-up factor for the evolution of the process. All of these are consequences of the correlations in the model, which are absent in the BTM. 

\subsubsection{Supercritical Liouville Brownian motion}
\label{sss:LBM}
The measure $\mu^{(\beta)}_N(N \cdot)$ may also be seen as a discrete version (up to  normalization) of the Liouville quantum gravity measure (LQGM, also known as Gaussian Multiplicative Chaos - GMC) associated to the continuum Gaussian free field (CGFF) on $V$. The latter is {\em formally} defined as a random measure whose Radon-Nykodym derivative with respect to Lebesgue measure is (up to formal normalization) $\rme^{-\beta h}$, where $h$ is the CGFF on $V$. In fact, one way of making sense of this formal definition is via a scaling limit similar to that in~\eqref{e:1.11} (see~\cite{rhodes2013gaussian}
for a general survey on LQGM/GMC and~\cite{barral2013gaussian} for a construction of general super-critical GMCs and their duality relation with the sub-critical ones). 

In a similar way, $X_N$ can be seen as the discrete analog of the Liouville Brownian motion (LBM) on $V$, which is formally a Brownian motion $B \equiv (B(t) :\: t \geq 0)$ time-changed by the inverse of the process $t \mapsto \int_{s=0}^t \rme^{-\beta h_{B(s)}} \rmd s$, where $h$ is again the CGFF on $V$. Existence of this object as a continuous strongly Markovian process has been shown in the subcritical case $\beta < \alpha$~\cite{B_lbm,GRV_lbm} and in the critical case $\beta = \alpha$ ~\cite{rhodes2015liouville}. Our work may be seen as the corresponding construction of LBM in the supercritical regime $\beta > \alpha$. It should be stressed that unlike in the case $\beta \leq \alpha$, our object is not constructed as a measurable function of the field $h$ and the motion $B$, but rather only as a distributional limit or, alternatively, by specifying its law directly. Both the LQGM/GMC and the LBM are fundamental objects in the physical theory of two-dimensional Liouville quantum gravity, whose mathematical formalization is an ongoing task.

\subsubsection{Particle in random media}
An alternative point of view is that of a particle whose motion is governed both by thermal activation and the disordered media in which it diffuses. In this context, it is perhaps more natural to study the process on an infinite domain (using the pinned DGFF on $\bbZ^2$ as the underlying media, for instance), where there is no stationary distribution. It is believed~\cite{carpentier2001glass, castillo2001freezing} that logarithmic correlations provide precisely the right balance between the the depth (energy) and number (entropy) of valleys in the environment $-h_{N}$, giving rise to phenomena not present when correlations decay slower, as in Sinai's walk on $\bbZ$~\cite{sinai1983limiting}, or faster, as in the random conductance model~\cite{biskup2011recent}. This model was recently studied in \cite{biskup2016return}, where predictions from \cite{carpentier2001glass,castillo2001freezing} were partially confirmed.

\subsection{Heuristics and outline}
\subsubsection{The extremal picture}
\label{ss:ExtremalPic}
The emergence of $Y^{(\beta)}$ as the scaling limit of $X_N$ becomes more clear if one looks at the extremal structure of the field for large $N$. Thanks to recent progress in this area, this structure is now well understood, and we proceed to describe the relevant results in this theory.

The principle extreme value, namely the global maximum of the field, was studied by Bramson, Ding and Zeitouni~\cite{bramson2016convergence}, who showed the existence of a random variable $M^*$, finite almost surely, such that with $m_N$ as in~\eqref{e:1.5},
\begin{equation}	
\label{eq:max_cvgs}
\max_{x \in V_N^*}\, (h_{N,x} - m_N) \Longrightarrow M^* \quad \text{as } N \to \infty \,.
\end{equation}
(This convergence of the centered maximum is the fruit of a long effort, dating back to a few decades before, and we invite the reader to consult~\cite{bramson2016convergence} for an historical overview.)

Other extreme values can be studied by defining for $A \subseteq \bbR$ and $N \geq 1$ the set
\begin{equation}
\label{eq:1.10}
\Gamma_N(A) := \big\{ x \in V_N :\: h_{N,x} - m_N \in A \big\} \,.
\end{equation}
Writing $\Gamma_N(v)$ as a shorthand for the extreme superlevel set $\Gamma_N([v, \infty))$, Ding and Zeitouni showed in~\cite{ding2014extreme} that for all $v \in \bbR$,
\begin{equation}
\label{eq:trap_separation}
\lim_{r\to \infty} \limsup_{N \to \infty}\, \bbP \Big( \exists x, y \in \Gamma_N\big(-v\big) : r < \|x-y\| < N /r \Big) = 0 \,.
\end{equation}
Thus, with high probability, the extreme values congregate in {\em clusters} of $O(1)$ diameter, which are  $N/O(1)$ apart.

This clustering of extreme values, which is a consequence of the short range correlations of the field, motivates the structured form of the extremal process $\eta_{N,r}$ in~\eqref{e:1.5}. This  process captures all extreme values by recording the location and height of the local maximum in each extremal cluster (the first two coordinates), and then separately the relative heights of all extreme values around it (the third coordinate).

In view of the convergence of $\eta_{N,r}$ to $\eta$ as given by $\eqref{eq:thm1}$, we see that the heights of the cluster maxima (centered by $m_N$) asymptotically form a Poisson point process with the exponential intensity $|Z| \rme^{-\alpha v} \rmd v$. The locations of these cluster maxima are asymptotically i.i.d. and chosen, after scaling by $N$, according to $\wh{Z}$. Lastly, relative to the cluster maximum, the field around each cluster is asymptotically chosen in an i.i.d. fashion according to the law $\nu$.

This asymptotic picture suggests that the size of the extreme superlevel sets satisfies $\log |\Gamma(-v)| \sim \alpha v$ as $v \to \infty$ with high probability for large $N$ (see Proposition~\ref{prop:tightness}). This indicates that when (and only when) $\beta > \alpha$, the  measure $\mu^{(\beta)}_N$ from~\eqref{e:1.11} concentrates on vertices corresponding to the extreme values of $h_N$. In view of the definition of both $\mu_N^{(\beta)}$ and $\eta_{N,r}$, it then holds that
\begin{equation}
\label{e:1.17}
\rme^{-\beta m_N} \mu_N^{(\beta)}(N \cdot) \approx \sum_{(x,v,\omega) \in \eta_{N,r}} \rme^{\beta v} \Big( \sum_{y \in \rmB_r} \rme^{-\beta \omega_y} \delta_{(x+y/N)} \Big) 
\approx 
\sum_{(x,v,\omega) \in \eta_{N,r}} \Big( \rme^{\beta v} \sum_{y \in \rmB_r} \rme^{-\beta \omega_y}  \Big) \delta_x \,. 
\end{equation}
Recalling that the mean holding time of $X_N$ at $x$ is $\mu_N^{(\beta)}(x)$ we see that each cluster $(x/N,v,\omega) \in \eta_{N,r}$ traps the walk for a time proportional in mean to $\rme^{\beta v} \sum_{y \in \rmB_r} \rme^{-\beta \omega_y}$.

Since $\eta_{N,r}$ converges to $\eta$, the right hand side in~\eqref{e:1.17} is approximately $\sum_{(x,t) \in \wt{\chi}^{(\beta)}} t \delta_x$, where $\wt{\chi}^{(\beta)}$ is obtained from $\eta$ via the pointwise transformation:
\begin{equation}
(x,v,\omega) \longmapsto \Big(x\,, \rme^{\beta v} \sum_{y \in \bbZ^2} \rme^{-\beta \omega_y} \Big) \,.
\end{equation}
But then an elementary calculation (Proposition~\ref{prop:depth_close}) using the finiteness of the integral in~\eqref{eq:kappa} shows that $\wt{\chi}^{(\beta)}$ has the law of $\chi^{(\beta)}$ described in~\eqref{eq:chi}, which is consistent with the first statement in~\eqref{e:1.11}. Thus the atoms of $\chi^{(\beta)}$ indeed encode the trapping landscape for $X_N$ in the limit.

As for the K-process, since the traps (or clusters) are $N/O(1)$ apart and $O(1)$ in diameter, standard random walk theory on the two-dimensional torus can be used to show that the number of returns to a trap, before exiting a ball of radius $O(N / \log N)$ around it, scales to an exponentially distributed random variable, and that following this exit, the next chosen trap is drawn approximately uniformly. This shows that for large $N$ the process $X_N$ exhibits the uniform trap hopping dynamics, which is characteristic of the spatial K-process. 

We remark that such random walk analysis has been done in the case of the BTM on the two-dimensional torus, both in in-equilibrium~\cite{JLT} and pre-equilibrium timescales~\cite{arous2008arcsine}. In this case, traps are single vertices (not finite clusters) and are scattered uniformly in $V_N^*$ (not according to $\wh{Z}$). Nevertheless, this analysis essentially still applies (with small modifications), and we therefore make use of some of the statements from these works, notably the uniform selection of traps (Lemma~\ref{lem:ord_close_pre}). 

\subsubsection{Outline of the paper}
Let us describe the structure of the rest of the paper. In Section~\ref{sec:CloseK}, we construct a {\em spatial pre K-process} from a given spatial K-process by forgetting all but the deepest $M$ traps. We then show that this process becomes close in law to the original process when $M$ is large. Section~\ref{s:DGFF} is devoted to characterizing the limiting trapping landscape of the DGFF, rigorizing the heuristic picture in the previous subsection. In Section~\ref{s:TraceClose}, we introduce the {\em trace process}. This process ``fast-forwards'' through vertices of $V_N^*$ which are not close to the $M$ deepest traps. As in Section~\ref{sec:CloseK}, the trace process is shown to be close in law to $X_N$ when $M$ is taken large. 

Next, in Section~\ref{s:trap_hop}, we study the trace process and show that it exhibits the uniform trap hopping dynamics of the spatial K-process, driven by the limiting trapping landscape of the DGFF. The outcome of this subsection is a coupling between the trace process and the spatial pre  K-process, in which they stay close to each other with high probability. In Section~\ref{s:final}, we combine this with the outcomes of Section~\ref{sec:CloseK} and~\ref{s:TraceClose}, and use standard weak convergence theory to complete the proof. Finally, Appendix~\ref{s:app} includes bounds on the discrete Green function in two dimensions both on the torus and in $\bbZ^2$.

\section{From the K-Process to the Pre K-Process}
\label{sec:CloseK}

It is convenient to introduce a simpler variant of the K-process, in which we keep only a finite number of the deepest traps. We call such a process a \emph{pre K-process} (often called a truncated K-process in the literature). The main result of this section is that the pre K-process is a good approximation to the full K-process, when the number of traps is taken to be large enough. We first consider a deterministic sequence of locations and depths and then treat the $\chi$-driven version of this process.

\subsection{Closeness of K-processes and pre K-processes}

Let $(\xi_k, \tau_k)_{k \geq 1}$ be a fixed collection of locations and depths, with the depths summable. Given $M \in \bbN$ and using the definitions and notation from the construction of the full K-process in Subsection~\ref{ss:SetupAndResults}, define the \emph{pre clock-process} $( T_M(u) : u \geq 0)$ as 
\begin{equation}
T_M \equiv T_M^{(\tau)}(u) := \sum_{k=1}^M \tau_k \sum_{j=1}^{A_k(u)} e_j^{(k)} \,,
\label{eq:2_pre_clock}
\end{equation}
and then the corresponding \emph{pre K-process} as
\begin{equation}
K_M(t) \! = \! K_M^{(\tau)}(t) \! := \! \left\{
	\begin{array}{ll}
	k	\!		& \text{if } t\! \in \! \big[T_M(\sigma_{j}^{(k)}\!-), T_M(\sigma_{j}^{(k)})\big) 
					\text{ for some } j \geq 1 \text{ and {\em one} } 1\!  \leq \! k \leq \! M , \\
	\infty	\!	& \text{otherwise.}
	\end{array} \right.
\label{eq:pre_K_def}
\end{equation}
Observe that $K_M(t) \neq \infty$ for all $t$ almost surely. In fact, $K_M$ is a random walk on the complete graph with $M$ vertices (self-loops included), with an exponentially distributed holding time with mean $\tau_k$ at the $k$-th vertex. From $K_M$ and the locations $\xi$, define the \emph{spatial pre K-process} $Y^{(\xi, \tau)}_M$:
\begin{align}
\label{eq:Ydef}
Y_M(t) = Y_M^{(\xi,\tau)}(t) := \xi_{K_M(t)}.
\end{align}

We start with a simple observation:
\begin{lem} Let $\eps, \mathfrak{t} >0$. There is $U_0( \eps, \mathfrak{t}, \tau) >0$ so that with probability at least $1- \eps$, for all $M \in \bbN$, and whenever $U \geq U_0$ we have $T_M(U) \geq \mathfrak{t}$. 
\label{lem:pre_K_interval_cover}
\end{lem}
\begin{proof}
This follows immediately from $T_M(u) \geq T_1(u) \to \infty$ as $u \to \infty$ which holds almost surely in light of, for instance, the strong law of large numbers.
\end{proof}

We now show that $M$ can be taken large enough so that the spatial K-process defined in~\eqref{eq:spatial_K} and the corresponding spatial pre K-process defined in~\eqref{eq:Ydef} agree for most of a given finite time interval. Below, we write $\Leb(A)$ for the Lebesgue measure of a Borel set $A \subset \bbR$. 

\begin{lem} Let $\eps, \mathfrak{t}>0$. There is $M_0(\epsilon,\mathfrak{t},\tau) \in \bbN$ so that with probability at least $1-\epsilon$, whenever $M \geq M_0$, we have $\Leb( B_M(\mft)) \leq \eps$, where \begin{align}
B_M(\mathfrak{t}) := \big\{ t \in [0,\mathfrak{t}] : Y^{(\xi, \tau)}(t) \neq Y^{(\xi, \tau)}_M(t) \big\}\,
\label{eq:K_pre_K_bad_times}
\end{align}
\label{lem:K_pre_K}
\end{lem}

\vspace{-10mm}

\begin{proof} It suffices to work with pre K-processes in place of spatial pre K-processes, and we do so throughout the proof. Let $M, M' \in \bbN$, to be chosen later with $M' < M$. For $\eps, \mathfrak{t} >0$, use Lemma~\ref{lem:pre_K_interval_cover} to extract $U_0$ depending on these parameters and the depths $\tau$, so that with probability at least $1-\eps$ we have $T_m(U_0) \geq \mathfrak{t}$ for all $m \in \bbN$. Work with this $U_0$ and within this high probability event, denoted $\mathcal{E}$, for the remainder of the proof. 

For $m \in \bbN$ let $S(m)$ be the tail of the full clock process at time $U_0$:
\begin{align}
S(m) := \sum_{k=m}^\infty \tau_k \sum_{j=1}^{A_k(U_0)} e_j^{(k)}, 
\end{align}
and let $J(m)$ be the number of jumps made by the first $m$ Poisson processes $A_k$ at time~$U_0$:
\begin{align}
J(m) := \sum_{k=1}^{m} A_k(U_0) \,.
\end{align}
For $\delta >0$, choose $M'$ large enough so that $\bbE ( S(M') ) / U_0 \equiv \sum_{k=M'}^\infty \tau_k \leq \delta^2$ and then choose $M$ so that $\bbE( S(M)) / U_0 \equiv \sum_{k=M}^\infty \tau_k \leq \delta^3 / M'$. Define the event 
\begin{align}
\mathcal{E}' := \big\{S(M') < \delta U_0\big\} \cap \big\{  J(M') < \delta^{-1} M' U_0 \big\}\cap \big\{ S(M) < \delta^2 U_0 / M' \big\} \,. 
\label{eq:K_pre_K_3}
\end{align}
By Markov's inequality, $\bbP( \mathcal{E}' ) \geq 1 - 3\delta$. Work also within $\mathcal{E}'$ for the remainder of the proof. 

Let $\{\sigma_i\}_{i=1}^\infty$ be the jump times of $K_M$, and let $\ell \in \bbN$ be smallest possible so that $\sigma_\ell  > \mathfrak{t}$. Setting $\sigma_0 = 0$, the intervals $\{ I_i \}_{i=1}^\ell \equiv \{ [\sigma_{i-1}, \sigma_i ) \}_{i=1}^\ell$ are disjoint and cover of $[0,\mathfrak{t}]$. The process $K_M$ is constant on each interval, and we let $K_M(I_i)$ denote the value of $K_M$ on the interval $I_i$. Call an interval $I_i$ \emph{extremely deep} (ED) if $K_M(I_i) \in \{1, \dots M'\}$, and call it \emph{moderately deep} (MD) otherwise.

If $[\sigma_{i-1}, \sigma_i)$ is extremely deep, then $K$ and $K_M$ agree on $[\sigma_{i-1}, \sigma_i) \setminus [ \sigma_{i-1}, \sigma_{i-1} + S(M) ]$, and thus,
\begin{align}
B_M(\mathfrak{t}) \subset \left( \bigcup_{ I_i \text{ is ED}} [\sigma_{i-1} , \sigma_{i-1} + S(M) ] \right) \cup \left( \bigcup_{ I_i \text{ is MD} } I_i \right)\,.
\end{align}
Within $\cal{E}$, the number of extremely deep intervals is at most $J(M')$, and the total length of the moderately deep intervals is at most $S(M')$. Thus,
\begin{align}
\Leb(B_M(\mathfrak{t}) ) \leq J(M') S(M) + S(M') \,.
\end{align}
By \eqref{eq:K_pre_K_3}, on $\mathcal{E}'$ we have $\Leb(B_M(\mathfrak{t})) \leq \delta U_0^2 + \delta U_0$. Taking $\delta = \epsilon / (1+U_0)^2$, we obtain $\Leb(B_M(\mathfrak{t}) ) \leq 2\epsilon$ with probability at least $1 - (3\delta + \epsilon) \geq 1 - 4\epsilon $. Taking $\epsilon/4$ for $\epsilon$ in the first place completes the proof.\end{proof}


\subsection{Closeness of $\chi$-driven processes}
For $M \geq 1$, we now define the $\chi$-driven spatial {\em pre} K-process, denoted $Y_M^{(\beta)} = \big(Y_M^{(\beta)}(t) : t \geq 0 \big)$. This process is defined conditionally on $\chi^{(\beta)}$, just as $Y^{(\beta)}$ was defined in Subsection~\ref{ss:SetupAndResults}, using the spatial pre K-process $Y_M^{(\xi^{(\beta)}, \tau^{(\beta)})}$ in place of the full one. 

We wish to show that a statement similar to the one in Lemma~\ref{lem:K_pre_K} holds also for the $\chi$-driven processes. We first establish that the $\chi$-driven spatial K-process is well-defined, a task postponed from Section~\ref{ss:SetupAndResults}. 

\begin{lem}
\label{lem:finite_tau}
Let $\beta > \alpha$ and let $\chi^{(\beta)}$ be defined as in~\eqref{eq:chi}. Then one may order the atoms of $\chi^{(\beta)}$ in descending order of their second coordinate. Moreover, if $\big((\xi^{(\beta)}_k, \tau^{(\beta)}_k)\big)_{k=1}^\infty$ denotes such an ordering, then $\sum_k \tau^{(\beta)}_k < \infty$ almost surely.
\end{lem}

\begin{proof}
Suppressing the superscript $\beta$, first note that the almost-sure finiteness of $Z$ and integrability of of $t^{-1-\alpha/\beta}$ on $[t, \infty)$ for all $t > 0$ implies $\chi(V^* \times [t, \infty)) < \infty$ for all $t > 0$, almost surely. This shows an ordering is indeed possible. To show summability, 
decompose the sum of the $\tau_k$ as
\begin{align}
\sum_{k \geq 1} \tau_k = \sum_{k\geq1} \tau_k \1_{(1,\infty]}(\tau_k) + \sum_{k\geq1} \tau_k \1_{[0,1]}(\tau_k) \,,
\label{eq:finite_tau_1}
\end{align}
and observe that by the argument above, the first sum has finitely many finite terms and hence must be finite. To handle the second sum on the right side of \eqref{eq:finite_tau_1}, take conditional expectation with respect to $Z$:
\begin{equation}
\bbE \Big(\sum_{k=1}^\infty \tau_k 1_{[0,1]}(\tau_k) \,\Big|\, Z \Big)
= \big( \kappa_\beta |Z| \big)^{\beta/ \alpha} \int_{0}^1 t^{-\alpha/\beta} \rmd t \,.
\end{equation}
This is finite almost surely since $\beta > \alpha$ and $|Z|<\infty$ with probability $1$.
\end{proof}

The next proposition establishes closeness of the $\chi$-driven processes, and is the main product of this section.
\begin{prop}
Fix $\beta > \alpha$ and let $\epsilon, \mathfrak{t} >0$ be given.
There is $M_0(\beta, \eps, \mathfrak{t}) \in \bbN$ so that with probability at least $1-\epsilon$, whenever $M \geq M_0$,
\begin{align}
\label{e:2.12}
\Big\| Y^{(\beta)}(\cdot) - Y_M^{(\beta)}(\cdot) \Big\|_{\rmL( [0,\mathfrak{t}],V^* )} < \epsilon \,,
\end{align}
\label{prp:K_pre_K}
\end{prop}

\begin{proof} Let $B_M^{(\beta)}$ be defined conditionally on $\chi^{(\beta)}$ as in \eqref{eq:K_pre_K_bad_times} with $(\xi, \tau) = ( \xi^{(\beta)}, \tau^{(\beta)})$. By Lemma \ref{lem:finite_tau}, we have $\sum_{k\geq1} \tau^{(\beta)}_{k} < \infty$ almost surely, and it follows from Lemma \ref{lem:finite_tau} and the bounded convergence theorem that $\Leb(B^{(\beta)}_M(\mathfrak{t}))$ tends to zero in probability as $M \to \infty$. It remains to observe that the left hand side of~\eqref{e:2.12} is bounded from above by $\Leb(B^{(\beta)}_M(\mathfrak{t}))$ times the diameter of $V^*$.
\end{proof}

\section{The Trapping Landscape: DGFF Extremes}\label{s:DGFF}

At in-equilibrium timescales, large $N$ and low temperatures, the dynamics are effectively governed by the extreme values of the underlying Gaussian free field, which determine the trapping landscape. In this section, we introduce the formal notion of a trap and collect various results concerning their joint geometry. The proofs of these results are mostly straightforward adaptations of corresponding statements concerning the extrema of $h_N$. 

\subsection{Defining the traps}
Let us introduce some notation. For $r > 0$ define the set of $r$-local maxima of $h_N$ by
\begin{equation}
\Lambda_{N}(r) := \left\{ x \in V_N :   h_{N,x} = \max_{y \in \rmB_r} h_{N,\,x+y} \right\} \,,
\label{eq:Lambda_def}
\end{equation}
where we recall that $\rmB_r = \rmB_r(0)$ is the open Euclidean ball of radius $r$ around $0$ in $\bbZ^2$.	
As $X_N$ becomes localized near each $x \in \Lambda_N(r)$, we shall refer to $x \in \Lambda_N(r)$ (and sometimes also to $\rmB_r(x)$) as a \emph{trap} with corresponding depth 
\begin{equation}
\label{eq:tau_def}
\tau_r(x) \equiv \tau^{(\beta)}_r(x) := \sum_{y \in \rmB_r} \rme^{\beta h_{N,\,x +y} } \,.
\end{equation}
Let $x_{N,k}$ denote the $k$-th trap in an enumeration of the elements of $\Lambda_N(r)$ in decreasing order of their depth. For $M \in \bbN$, define the set of \emph{deep traps} as
\begin{equation}
\Lambda_{N}(r,M) \equiv \Lambda^{(\beta)}_{N}(r,M) := \big\{x_{N,1}, \dots, x_{N,M} \big\}\,.
\label{eq:Lambda_M_def}
\end{equation}
Thus, $\tau_r(x_{N,1}) > \tau_r(x_{N,2}) > \dots > \tau_r(x_{N,M}) > \max_{y \in \Lambda_N(r) \setminus \Lambda_N(r,M)} \tau_r(y)$. For convenience, we set $\ol{\Lambda}_N(r,M) := \rmB_r(\Lambda_N(r,M))$, where in general we write $\rmB_r(A)$ for $\cup_{x \in A} \rmB_r(x)$. The remaining vertices are denoted as $\ol{\Lambda}^\rmc_{N}(r,M) := V_N \setminus \ol{\Lambda}_{N}(r,M)$. 

\subsection{Dominance of deep traps}
Below, in Proposition~\ref{prop:tightness}, we show that the Gibbs distribution in \eqref{eq:1.9} is asymptotically ``carried'' by deep traps when $\beta > \alpha$.

\begin{prop} Let $\beta > \alpha$. Then for any $\eps >0$, there is $M_0(\eps,\beta) \in \bbN$, $r_0(M, \eps, \beta)$ and $N_0(M, r, \eps, \beta) \in \bbN$ so that for all $M \geq M_0$, $r > r_0$ and $N \geq N_0$ with probability at least $1-\eps$, 
\begin{align}
\label{e:3.5}
\sum_{x \in \ol{\Lambda}_N^\rmc (r,M)} \rme^{\beta (h_{N,x} - m_N)} < \eps \,.
\end{align}
\label{prop:tightness}
\end{prop}

We prove two lemmas to deduce Proposition~\ref{prop:tightness}.  Recall the notation $\Gamma_N(A)$ defined in \eqref{eq:1.10}, and recall also the shorthand $\Gamma_N(v)$ for $\Gamma_N( [v,\infty))$; let us also abbreviate $V_N \setminus \Gamma_N(v)$ as $\Gamma_N^\rmc(v)$. 

\begin{lem}
\label{l:3.2}
For all $\epsilon > 0$, there is $u_0(\eps) >0$ and $N_0(\eps, u)$ so that $u \geq u_0$ and $N \geq N_0$ imply that with probability at least $1-\eps$, 
\begin{align}
\sum_{x \in \Gamma_N^\rmc(-u)} \rme^{\beta ( h_{N,x} - m_N )} < \eps 
\label{e:3.6}
\end{align}
\end{lem}
\begin{proof}
By Proposition 6.8 in~\cite{biskup2016full}, there is $\beta' \in (\alpha, \beta)$ and $c > 0$ so that for all $v >0 $ and $N \in \bbN$ large enough,
\begin{equation}
\label{e:3.7}
\bbP  \big( \big|\Gamma_N(-v)\big| > \rme^{\beta' v} \big) \leq \rme^{-cv} \,.
\end{equation} 
Summing the above probabilities along $v = k \in \bbN$ and using the union bound, we see that for any $\epsilon > 0$ there exists $k_0 \in \bbN$ so that for all $N$ large enough, with probability at least $1-\epsilon$, we shall have $|\Gamma_N(-k)\big| \leq \rme^{\beta' k}$ for all $k \geq k_0$. But then, if $\lfloor u \rfloor \geq k_0$, the sum in \eqref{e:3.6} is bounded above by
\begin{equation}
\sum_{k \geq \lfloor u \rfloor} \big| \Gamma_N(-(k+1)) \setminus \Gamma_N(-k) \big| \rme^{-\beta k} \leq \sum_{k \geq \lfloor u \rfloor} \big| \Gamma_N(-(k+1)) \big| \rme^{-\beta k}
\leq C \rme^{-(\beta - \beta') u} \,,
\end{equation}
which will be smaller than $\epsilon$ for all $u$ large enough.
\end{proof}

The next result relates deep traps to the extreme superlevel sets of the field. 

\begin{lem}
\label{l:3.3}
Let $u >0$ and let $\eps >0$. There are $M_0(\eps, u)$, $r_0(\eps,M,u)$ and $N_0(\eps, M, r, u)$ so that $M \geq M_0$, $r \geq r_0$ and $N \geq N_0$ imply 
\begin{equation}
\label{e:3.9}
\bbP\big( \Gamma_N(-u)  \subseteq \ol{\Lambda}_N(r,M) \big) \geq 1 -\eps \,.
\end{equation}
\end{lem}

\begin{proof} Fix $u >0$. For $M \in \bbN$, $r >0$ and $v \geq u$, let $\cA_1$ be the event that~\eqref{e:3.6} holds with $v$ in place of $u$ and with $\epsilon = \rme^{-\beta u}$, namely:
\begin{align}
\cA_1 := \Big\{ \textstyle\sum_{x \in \Gamma_N^\rmc(-v)} \rme^{\beta ( h_{N,x} - m_N )} < \rme^{-\beta u} \Big\}\,,
\label{eq:calA1}
\end{align}
 and define events $\cA_2$ and $\cA_3$ as follows:
\begin{align}
\cA_2 & := \Big\{ \big|\Gamma_N(-v)\big| \leq M \Big \} \,, \\
\cA_3 & := \Big\{ \forall x,y \in V_N , \: x,y \in \Gamma_N(-u)  \text{ implies } \|x-y\| 
	\notin (r, N/r) \Big\} \,.
\end{align}
Observe that the event in \eqref{e:3.9} contains the intersection $\cA_1 \cap \cA_2 \cap \cA_3$. Hence, by the union bound it is enough to show that the probability of the complement of each goes to $0$ in the limit when $N \to \infty$ followed by $r \to \infty$ then $M \to \infty$ and finally $v \to \infty$. Indeed, this holds for $\cA_1^\rmc$ in light of Lemma~\ref{l:3.2}, for $\cA_2^\rmc$ by Proposition 6.8 in~\cite{biskup2016full} as in~\eqref{e:3.7} and finally for $\cA_3^\rmc$ thanks to~\eqref{eq:trap_separation}. 
\end{proof}




Combining these lemmas, we can easily present

\begin{proof}[Proof of Proposition~\ref{prop:tightness}]
Fix $\beta > \alpha$ and $\epsilon > 0$ and use Lemma~\ref{l:3.2} to find $u > 0$ such that the event in~\eqref{e:3.6} occurs with probability is at least $1-\epsilon/2$ for all $N$ large enough. Then using Lemma~\ref{l:3.3} find $M_0$, $r_0(M)$ and $N_0(r,M)$ such that whenever $M \geq M_0$, $r \geq r_0(M)$ and $N \geq N_0(r,M)$ the event in~\eqref{e:3.9}
occurs with probability at least $1-\epsilon/2$. But on the intersection of the last two events, which has probability at least $1-\epsilon$, the inequality in~\eqref{e:3.5} holds.
\end{proof}

\subsection{Separation of deep traps}\label{sec:3_sep}

Next we address the separation of deep traps. Let $\partial V_N$ denote the vertices of $\bbZ^2$ which are neighbors of vertices in $V_N$, but are not in $V_N$, and define $r_N := N / \log N$. We say the field $h_N$ is \emph{$(r,M)$-separated} if all deep traps are at least a distance of $r_N$ from one another and from the boundary of $V_N$, that is:
\begin{align}
\min_{\substack{ \text{$x,y \in \Lambda_N(r,M)$}\\\text{$x \neq y$}}} \| x - y \| \geq r_N \hspace{5mm} \text{ and } \  \min_{x \in \Lambda_N(r,M),\, z \in \pa V_N} \| x - z \| \geq r_N \,.
\label{eq:well-separated}
\end{align}
Henceforth we denote the event that $h_N$ is $(r,M)$-separated by $\mathcal{S}_N(r,M) \equiv \mathcal{S}_N^{(\beta)} (r,M)$.

\begin{prop} Fix $\beta > \alpha$ and let $\e >0$ and $M \in \bbN$ be given. There is $r_0 = r_0(M, \eps, \beta)$ and $N_0(r) = N_0(r, M, \eps, \beta)$ so that $r \geq r_0$ and $N \geq N_0$ together imply that $h_N$ is $(r,M)$-separated with probability at least $1-\eps$.
\label{prop:extreme_sep}
\end{prop}

We first show that when $r$ and $M$ are fixed, we may choose $u$ large enough so that $\Lambda_N(r,M) \subseteq \Gamma_N([-u,u])$ with high probability. 

\begin{lem}
\label{l:3.5}
For all $M \in \bbN$ and $r >0$. There is $u_0(M,r) >0$ and $N_0(M,r,u) \in \bbN$ so that when $u \geq u_0$ and $N \geq N_0$,
\begin{align}
\label{e:3.14}
\bbP( \Lambda_N(r,M) \subseteq \Gamma_N( [-u,u]) ) \geq 1- \eps 
\end{align}
\end{lem}

\begin{proof} Fix the parameters $M$ and $r$. Let $r' \geq r$ and choose $v \in [0,u]$. Let $\cA_1$ be the event from~\eqref{eq:calA1}. Recall the definition of $\eta_{N,r}$ from~\eqref{e:1.5}, and define the following events:
\begin{align}
\cA_4 & := \Big\{ \eta_{N,r'}\big([0,1]^2 \times [-v, 0] \times [0,\infty)^{\bbZ^2} \big) \geq M \Big\} \,, \\
\cA_5 & := \Big\{ \Gamma_N(u) = \emptyset \Big\} \,.
\end{align}
Observe that for any $v \in [0,u]$, the event in \eqref{e:3.14} contains the intersection $\cA_1 \cap \cA_4 \cap \cA_5$. As before, it is therefore sufficient to show that the complementary probabilities go to zero when $N \to \infty$, followed by $u \to \infty$, then $r' \to \infty$ and finally $v \to \infty$. 

Indeed, this is true for $\cA_1$ thanks to Lemma~\ref{l:3.2}. For $\cA_4$ this holds in light of~\eqref{eq:thm1} which implies that $\bbP(\cA_4^\rmc)$ tends, when $N \to \infty$ followed by $r' \to \infty$, to the probability that a Poisson random variable with conditional mean $Z\alpha^{-1} (\rme^{\alpha v} - 1)$ is smaller than $M$. In light of the almost-sure finiteness of $Z$, this probability will go to $0$ when $v \to \infty$. Finally, $\cA_5$ is handled by~\eqref{eq:max_cvgs}. 
\end{proof}

It is now a short step to 
\begin{proof}[Proof of Proposition~\ref{prop:extreme_sep}]
Fix $\beta > \alpha$ and let $\epsilon > 0$ and $M > 0$ be given. Use Lemma~\ref{l:3.5} to find $u>0$ so that the event in~\eqref{e:3.14} holds with probability at least $1-\epsilon / 3$ for all $N$ large enough. By~\eqref{eq:trap_separation} we may then find $r$ large enough so that for all $N \geq N_0(r)$, any $x,y \in \Gamma_N(u)$ not within Euclidean distance $r$ from each other, must be at least $N/r$ from one another with probability at least $1-\epsilon/3$. Thus, any distinct $x, y \in \Lambda_N(r) \cap \Gamma_N(u)$ must be at least $r_N$ apart for large $N$.

At the same time, recalling the law of $\eta$ from~\eqref{eq:thm1} and the fact that $Z$ does not charge the boundary of $[0,1]^2$ almost surely, we may find $\delta > 0$ so that with with probability at least $1-\epsilon/6$ there is no mass in $\big([0,1]^2 \setminus [\delta, 1-\delta]^2\big) \times [-u, u] \times [0,\infty)^{\bbZ^2}$ under $\eta$. In light of the convergence in~\eqref{eq:thm1}, this translates into the assertion that for all $r$ large enough and then $N$ large enough, $\Lambda_N(r) \cap ( V_N \setminus [r_N, N-r_N]^2) \cap \big(\Gamma_N(u) \setminus \Gamma_N(-u)\big)  = \emptyset$ with probability at least $1-\epsilon /3$. 

Combining the above, we see that $S_N(r,M)$ occurs with probability at least $1-\epsilon$, once $r$ and then $N$ are chosen sufficiently large.
\end{proof}

\subsection{Law asymptotics for the structure of deep traps}
The last task of this section is to address the convergence of the locations and depths of deep traps to the atoms of the process $\chi^{(\beta)}$ from~\eqref{eq:chi} in the limit of large $N$.

\begin{prop}
Fix $\beta > \alpha$. Then for all $M \geq 1$, as $N \to \infty$ followed by $r \to \infty$, we have
\begin{align}
\label{e:3.17}
\Big( x_{N,j} / N ,\, \rme^{-\beta m_N} \tau_r(x_{N,j}) \Big)_{j=1}^M 
\Longrightarrow  (\xi_j, \tau_j)_{j=1}^M \,,
\end{align}
where $( \xi_j, \tau_j)_{j=1}^M$ enumerate the first $M$ atoms of $\chi^{(\beta)}$, ordered in a decreasing manner according to the second coordinate. 
\label{prop:depth_close}
\end{prop}

\begin{proof} 
For any $R, r_0 \in (0, \infty)$, let us define the function $f_{R,r_0}$ from $[0,1]^2 \times (-\infty, R] \times [0, \infty)^{\rmB_{r_0}}$ to $[0,1]^2 \times (0, \infty)$ via 
\begin{equation}
f_{R,r_0}(x, h,\omega) := \Big(x,\, \sum_{y \in \rmB_{r_0}} \rme^{\beta (h - \omega_y)}
 \Big) \,.
\end{equation}
Notice that $f_{R,r_0}$ is continuous in the product topology and that preimages of compact sets under $f_{R,r_0}$ are compact. In particular, $f_{R,r_0}$ induces a continuous function $F_{R,r_0}$ from the space of Radon point measures on $[0,1]^2 \times (-\infty, R] \times [0, \infty)^{\rmB_{r_0}}$ 
to the space of Radon point measures on $[0,1]^2 \times (0, \infty)$ via
\begin{equation}
F_{R,r_0}(\eta) := \sum_{(x,h,\omega) \in \eta} \delta_{f_{R,r_0}(x,h,\omega)} \,,
\end{equation}
where continuity is with respect to the vague topology.

Let $\eta_{N,r}^{R,r_0}$ and $\eta^{R,r_0}$ denote, respectively, the restriction and proper projection of $\eta_{N,r}$ and $\eta$ from~\eqref{eq:thm1} onto $[0,1]^2 \times (-\infty, R] \times [0, \infty)^{\rmB_{r_0}}$. It follows from~\eqref{eq:thm1} and the fact that $[0,1]^2 \times (-\infty, R] \times [0, \infty)^{\bbZ^2}$ is stochastically continuous under $\eta$, that $\eta_{N,r}^{R,r_0}$ tends weakly to $\eta^{R,r_0}$ when $N \to \infty$ followed by $r \to \infty$. But then, continuity of $F_{R,r_0}$ implies that when $N \to \infty$ and then $r \to \infty$,
\begin{equation}
\label{e:3.28}
\chi_{N,r}^{R,r_0} := F_{R,r_0}(\eta_{N,r}^{R,r_0}) \,\,\Longrightarrow \,\,  F_{R,r_0}(\eta^{R,r_0}) =: \chi^{R,r_0}\,.
\end{equation}

Next we wish to take $r_0$ and $R$ to infinity in both $\chi_{N,r}^{R,r_0}$ and $\chi^{R,r_0}$. We begin with the former object; let $\epsilon > 0$ and use Lemma~\ref{l:3.2} to find $v$ so that for all large $N$, and with probability at least $1-\epsilon$,
\begin{equation}
\sum_{x \in \Gamma^\rmc_N(v)} \rme^{\beta (h_{N,x} - m_N)} < 
\epsilon \,.
\end{equation} 
Assuming that this event holds, we have that $\chi_{N,r}^{r_0}$ does not charge $N^{-1} \Gamma^\rmc_N(v) \times [\epsilon, \infty)$ for all $r_0 \leq r$. Furthermore, by~\eqref{eq:trap_separation}, we may find $r_0$ large enough, such that all $x,y \in \Gamma_N(v)$ are either at most $r_0$ or least $N/r_0$ apart with probability at least $1-\epsilon$, provided $N$ is large enough. Assuming that this event holds as well, we have for all $r \geq r_0$,
\begin{equation}
\sum_{x \in \Lambda_N(r) \cap \Gamma_N(v)}
\sum_{y \in \rmB_r \setminus \rmB_{r_0}} \rme^{\beta (h_{N,x + y} - m_N)}
\leq \epsilon \,,
\end{equation}
as all $x+y$ in the above sum are distinct and included in $\Gamma^\rmc_N(v)$.

This implies that if $\varphi \in C^\infty_0\big([0,1]^2 \times (\epsilon,R]\big)$, then for all $r_0$ and then $N$ large enough, we shall have 
\begin{equation}
\Big| \int \varphi(x, \tau)\, \chi_{N,r}^{R,r} (\rmd x \rmd \tau) - \int \varphi(x, \tau)\, \chi_{N,r}^{R,r_0 \wedge r} (\rmd x \rmd \tau) \Big| \leq \big\|\nabla \varphi\big\|_\infty\, \epsilon \,,
\end{equation}
with probability at least $1-2\epsilon$, where the norm on the right is the usual supremum norm.
But as $\epsilon$ and $\varphi$ were arbitrary, this shows that $\chi_{N,r}^{R,r_0 \wedge r}$ converges to $\chi_{N,r}^{R,r}$ vaguely in probability as $r_0 \to \infty$, uniformly in $N$ and $r$.

For a similar statement when $R \to \infty$, first use the almost-sure boundedness of the support of $\eta_{N,r}$ to define
\begin{equation}
 \chi_{N,r} := \lim_{R \to \infty} \chi_{N,r}^{R,r} = 
 	\sum_{x \in \Lambda_N(r)} \delta_{(x/N,\, \exp(-\beta m_N) \tau_r(x))} \,.
 \end{equation}
Using the tightness of the centered maximum, as implied by~\eqref{eq:max_cvgs}, for any $\epsilon > 0$ if $R$ is large then $\eta_{N,r}\big([0,1]^2 \times (R, \infty) \times [0,\infty)^{\rmB_r}\big) = 0$  with probability at least $1-\epsilon$ for all $r$ and $N$. On this event we must have $\chi_{N,r}^{R,r} = \chi_{N,r}$, and it follows that $\chi_{N,r}^{R,r}$ converges vaguely in probability to $\chi_{N,r}$ as $R \to \infty$, with the convergence uniform in $N$ and $r$. Altogether we have shown that when $r_0 \to \infty$ and then $R \to \infty$,
\begin{equation}
\label{e:3.33}
\chi_{N,r}^{R,r_0 \wedge r} \longrightarrow \chi_{N,r} \,
\end{equation}
vaguely in probability uniformly in $N$ and $r$.

Turning to $\chi^{R,r_0}$, recall that $\eta^{R,r_0}$ is, conditionally on $Z$, a Poisson point process on $[0,1]^2 \times (-\infty, R] \times [0, \infty)^{\rmB_{r_0}}$ and 
let $I^{R, r_0}_Z(\rmd x\,\rmd h\,\rmd \phi) := Z(\rmd x)\otimes \rme^{-\alpha h} \rmd h \otimes \nu(\rmd \phi)$ be its (conditional) intensity measure. An elementary ``change-of-variable'' calculation shows that 
\begin{equation}
\label{eq:3.32}
\big(I^{R,r_0}_Z f^{-1}_{R, r_0}\big) \big(\rmd x \, \rmd \tau \big) = 
Z(\rmd x) \otimes \tau^{-\alpha/\beta - 1} \kappa^{R,r_0}_{\beta}(\tau)  \rmd \tau \,,
\end{equation}
where $\tau \mapsto \kappa_\beta^{R, r_0}(\tau)$ is a positive function on $[0,\infty)$ which is pointwise increasing in $r_0$ and $R$, and which pointwise converges to $\kappa_\beta$ from~\eqref{eq:kappa}, when $r_0$ and $R$ tend to infinity.

It follows from the monotone convergence theorem that, conditional on $Z$, that the measure $\big(I^{R,r_0}_Z f^{-1}_{R, r_0}\big)(\rmd x\, \rmd \tau)$ tends vaguely to $Z(\rmd x) \otimes \kappa_\beta \tau^{-\alpha/\beta - 1} \rmd \tau$ when $r_0 \to \infty$ and $R \to \infty$. Since these are the conditional intensity measures of $\chi^{R, r_0}$ and $\chi^{(\beta)}$ respectively, which are both conditionally Poisson, it follows that when $r_0 \to \infty$ and $R \to \infty$,
\begin{equation}
\label{e:3.35}
\chi^{R, r_0} \Longrightarrow \chi^{(\beta)} \,.
\end{equation}
This is true conditionally on $Z$ almost surely, but then thanks to the bounded convergence theorem, also without the conditioning.

Combining~\eqref{e:3.28},~\eqref{e:3.33} and~\eqref{e:3.35} and 
using standard weak convergence (or metric space) theory, we see that $\chi_{N,r} \Longrightarrow \chi^{(\beta)}$ as $N \to \infty$ followed by $r \to \infty$. 
It remains to observe that $\big(x_{N,j}/N,\, \rme^{-\beta m_N} \tau_r(x_{N,j}) : j=1, \dots, M \big)$ as defined in~\eqref{eq:Lambda_M_def} are exactly the top $M$ atoms of $\chi_{N,r}$ ordered by their second coordinate. The convergence in~\eqref{e:3.17} will therefore follow by standard theory of point processes, provided we show that $\rme^{-\beta m_N} \tau_r(x_{N,1})$ is tight in $N$ and $r$. But since $\tau_r(x_{N,1}) \leq \sum_{x \in V_N} \rme^{\beta h_{N,x}}$, this follows from Theorem 2.6 in~\cite{biskup2016full}, which shows that $\sum_{x \in V_N} \rme^{\beta (h_{N,x} - m_N)}$ admits a proper limit in law when $N \to \infty$.
\end{proof}

\section{From the Random Walk to its Trace}
\label{s:TraceClose}

In this section, we introduce the {\em random walk trace process}, or \emph{trace process} for short. Conditionally on $h_N$, this process evolves like the original random walk $X_N$, except that it spends {\em zero} time at vertices not belonging to a deep trap. As with the pre K-process, the key output of this section is that the trace process can be made to be arbitrarily close to the original random walk in the $\|\cdot\|_{\rmL([0,\mft], V^*)}$-metric, provided that we consider enough deep traps. 

Let us now define this process explicitly. Fix $r > 0$ and $M \in \bbN$, and let the time spent by $X_N$ in $\ol{\Lambda}_N(r,M)$ up to time $t$ be defined via
\begin{align}
F_N^{(r,M)}(t) := \int_0^t \1_{\left\{ X_N(s) \in \ol{\Lambda}_N(r,M) \right\}} ds\,.
\end{align}
The \emph{random walk trace process} $X_N^{(r,M)}$ is constructed from $X_N$ via the time change: 
\begin{equation}
X_N^{(r,M)}(t) := X_N \Big( \overleftarrow{F}_N^{(r,M)}(t) \Big) \,,
\end{equation}
where $\overleftarrow{F}_N^{(r,M)}$ is generalized inverse of $F_N^{(r,M)}$, namely 
\begin{align}
\overleftarrow{F}_N^{(r,M)}(t) := \inf \big\{ s \geq 0 : F_N^{(r,M)}(s) \geq t \big\}.
\label{eq:4.2.5}
\end{align} 
Recall that in Theorem~\ref{thm:A}, the process $X_N$ is observed for time
$s_N \mft$ where,
\begin{equation}
s_N := g N^{2\sqrt{g} \beta}(\log N)^{1-3\sqrt{g}\beta/4} 
 \equiv g \rme^{\beta m_N} \log N \,.
\end{equation}
\begin{prop}
\label{prop:TraceClose}
 Let $\beta > \alpha, \eps >0$ and $\mathfrak{t} >0$. There are $M_0( \beta, \eps, \mathfrak{t} ) \in \bbN$, $r_0 (\beta, \eps, M, \mathfrak{t}) >0$ and $N_0(\beta, \eps, M, r, \mathfrak{t}) \in \bbN$ so that when $M \geq M_0$, $r > r_0$ and $N \geq N_0$, the event 
\begin{align}
\left\| \tfrac{1}{N} X_N(s_N \cdot) - \tfrac{1}{N} X_N^{(r,M)} (s_N \cdot) \right\|_{\rmL ( [0,\mathfrak{t}], V^* ) } < \epsilon
\end{align}
occurs with probability at least $1- \eps$.
\end{prop}
The proof of this proposition makes use of two constructions, that we introduce in Subsections~\ref{ss:RWClock} and~\ref{sec:macroscopic_jumps}. These constructions will also be used in other parts of the paper and we refer to them when necessary. Subsection~\ref{ss:TraceClose} includes the remaining argument needed to complete the proof of the proposition.

\subsection{From time to steps: the random walk clock process}
\label{ss:RWClock}

The process $X_N$ is a continuous time symmetric random walk with exponentially distributed holding times having mean $\rme^{\beta h_{N,x}}$ at vertex $x$. As such, we may construct $X_N$ using a discrete time simple random walk $\sfX_N = (\sfX_N(j) : j \geq 0)$ on $V_N^*$ and a collection of i.i.d. exponential random variables $\big(E_{N,j} : j \geq 1\big)$ independent of $\sfX_N$ and the field. Using these we can express $X_N$ as $X_N(s) = \sfX_N( ( \overset{{}_{\leftarrow}}{t}(s) -1 )^+)$, where
\begin{align}
t(n) &:= \sum_{j=0}^{n-1} \exp\left( \beta h_{N, \sfX_N(j)} \right) E_{N,j} \,,
\end{align}
and where $ \overset{{}_{\leftarrow}}{t}(s)$ is the generalized inverse of $t$ as in \eqref{eq:4.2.5}. The process $t(\cdot)$ is the {\em clock process} associated with $X_N$. It is no coincidence that the same name is also used in the literature for the process $T^{(\tau)}$ from~\eqref{eq:clock}, as both processes essentially serve the same purpose. Nevertheless, to avoid confusion, we refer to $t(\cdot)$ as the {\em random walk clock process}. 

The corresponding {\em trace clock process}, namely the clock process associated with $X_N^{(r,M)}$, is $\big(t^{(r,M)}(n) : n \geq 0\big)$ where,
\begin{equation}
t^{(r,M)}(n) := \sum_{j=0}^{n-1} \exp\left( \beta h_{N, \sfX_N(j)} \right) E_{N,j} \1
_{\left\{\sfX_N(j) \in \ol{\Lambda}_N(r,M) \right\}}  \,,
\label{eq:discrete_trace_clock_process}
\end{equation}
so that $X_N^{(r,M)}(s) = \sfX_N(( \overset{{}_{\leftarrow}}{t}{}^{(r,M)}(s) -1)^+)$ for $s \geq 0$. 

The advantage of using $\sfX_N$ and $t(\cdot)$ to describe $X_N$ is that the motion is decoupled from elapsed time. As such, ``spatial'' statements can be proved using simple random walk theory for $\sfX_N$, provided that we translate the elapsed time for $X_N$ to the number of steps of $\sfX_N$. The next lemma shows the corresponding number of steps of $\sfX_N$ is 
\begin{align}
\vart_N(n) := n \left\lceil N^2 \log N \right\rceil \,,
\end{align}
with $n$ properly chosen. In fact, we may choose $n$ so that the desired statement holds also for the $(r,M)$-trace process, and even uniformly in $r$ and $M$. 
\begin{lem} Let $\beta > \alpha$, $\eps >0$ and $\mathfrak{t} >0$. There is $n_0(\beta, \eps, \mathfrak{t}) \in \bbN$ so that for all $r >0$, $M \in \bbN$, all $N$ sufficiently large and $n \geq n_0$,
\begin{align}
t^{(r,M)} ( \vart_N(n)) \geq s_N \mft
\end{align}
occurs with probability at least $1- \eps$. 
\label{lem:2_elapsed_cover}
\end{lem}

\begin{proof} Let $x_N^*$ denote the argmax of $h_N$, and write $\vart_N(n)$ as $\vart$. Define
\begin{align}
\ell_\vart(x_N^*) = \sum_{k=0}^\vart \textbf{1}_{\{\sfX_N(k) = x_N^*\}} \,,
\end{align}
and for $(E_k)_{k \geq 1}$ a collection of i.i.d. mean $1$ exponentially distributed random variables independent of the $h_N$ and $\sfX_N$, also define the random variable
\begin{align}
t^*(\vart) := \exp(\beta h_{N,x_N^*}) \sum_{k=1}^{\ell_{\vart}(x_N^*) } E_k\,.
\end{align}
Letting $\sfG^N_n(0,y)$ be the expected number of visits $\sfX_N$ makes to $y$ in $n$ steps starting from $0$, by the ergodic theorem for $\sfX_N$ and the weak law of large numbers for the sequence $(E_k)_{k \geq 1}$, there is $N_0$ large so that with probability at least $1- \eps$, we have
\begin{align}
\frac{\ell_{\vart}( x_N^*)}{\sfG_{\vart}^N(0, x_N^*)} \geq 1/2\,,\hspace{10mm} \frac{1}{\ell_{\vart}( x_N^*)} \sum_{k=1}^{\ell_{\vart}( x_N^*)} E_k \geq 1/2 \,.
\label{eq:elapsed_cover_1}
\end{align}
when $N \geq N_0$. Invoke Lemma \ref{lem:2_green_estimate_one}, taking $N$ larger if necessary so that
\begin{align}
\frac{\sfG_{\vart}^N(0, x_N^*)}{\log N} \geq \frac{cn}{2} \,,
\label{eq:elapsed_cover_2}
\end{align}
and use \eqref{eq:elapsed_cover_1}, \eqref{eq:elapsed_cover_2} to obtain
\begin{align}
\bbP \Big( s_N^{-1} t^*(\vart) \geq \exp(\beta (h_{N,x_N^*} -m_N)) \frac{cn}{8g} \Big| h_N \Big) \geq 1- \eps \,.
\end{align}
Using the tightness of the centered maximum, per~\eqref{eq:max_cvgs}, choose $t >0$ according to $\eps$ so that
\begin{align}
\bbP \Big( s_N^{-1} t^*(\vart) \geq \frac{cn}{8g} \rme^{-\beta t} \Big) \geq 1-2\eps \,.
\end{align}
As $t^{(r,M)}( \vart)$ stochastically dominates $t^*(\vart)$, the proof is complete.
\end{proof}

\subsection{From steps to trap visits: macroscopic jumps}\label{sec:macroscopic_jumps}
To understand how $X_N$ and its trace $X_N^{(r,M)}$ move within the landscape, we further decompose the trajectory of $\sfX_N$ according to the visits it pays to deep traps. In light of Proposition~\ref{prop:extreme_sep}, we recall that such traps tend to be at least $r_N$ apart -- thus we think of a segment of the trajectory of $\sfX_N$ starting at the first entrance to $\rmB_r(x_{N,k})$ and ending at the first exit from $\rmB_{r_N}(x_{N,k})$ as a visit to the deep trap $x_{N,k}$. Between successive visits to deep traps, we say $\sfX_N$ makes a macroscopic jump.

To formalize the above construction, introduce the following sequence of stopping times, which are almost surely finite by recurrence of $\sfX_N$: 
\begin{align}
\sfR_1 &:= \inf \big\{ n \geq 0 : \sfX_N(n) \in \ol{\Lambda}_N(r,M)) \big\}  \\
\sfS_1 &:= \inf \big\{ n \geq \sfR_1 : \sfX_N(n) \notin \rmB_{r_N}\big(\Lambda_N(r,M)\big) \big\}
\end{align}
and for $k \geq 2$,
\begin{align}
\sfR_k &:= \inf \big\{ n \geq \sfS_{k-1} : \sfX_N(n) \in \ol{\Lambda}_N(r,M) \big\}  \\
\sfS_k &:= \inf \big\{ n \geq \sfR_k : \sfX_N(n) \notin \rmB_{r_N}\big(\Lambda_N(r,M)\big) \big\} \,.
\label{eq:stopping_times}
\end{align}
At each time $\sfR_k$, say $\sfX_N$ has made a \emph{$(r,M)$-macroscopic jump}. For $\vart \in \bbN$, define the number of $(r,M)$-macroscopic jumps made by $\sfX_N$ in $\vart$ steps as
\begin{align}
\sfJ_N^{(r,M)}(\vart) := \max \big\{ k\geq 0 \:: \sfR_k \leq \vart \big\} \,.
\label{eq:macro_jump_def}
\end{align} 

Next, we wish to show that within $\vart_N(n)$  steps, the number of macroscopic jumps is bounded from above with high probability uniformly in $N$. To do this, we need the next lemma, which is a nearly immediate consequence of (6.18) from \cite{JLT}. For what follows let $\sfH_N(A)$ be the hitting time of $A \subseteq V_N^*$ by $\sfX_N$ and write $\sfH_N(x)$ for $\sfH_N(\{x\})$. When we make general statements concerning $\sfX_N$, we will allow it to start from any vertex $x$, and use $\sfP_x^N$ and $\sfE_x^N$ to denote the underlying probability measure and expectation in this case.

\begin{lem} 
For all $\epsilon > 0$ and $r >0$, there exists $n_0 = n_0(\epsilon)$ and $N_0 = N_0(\epsilon , r)$, such that for all $N \geq N_0$,
\begin{align}
\sup_{{\substack{x,y\, \in\, V_N^*\\ d_N^*(x,y)\, >\, r_N / 2}}}\sfP_x^N \Big( \sfH_N( \rmB_r(y)) < \vart_N(n_0) \Big) \leq \epsilon \,.
\label{eq:A.1_jump_control_true}
\end{align}
\label{lem:jump_control_true}
\end{lem}

\begin{proof} 
It follows from recurrence of the simple random walk on $\bbZ^2$ that for any $r > 0$ there is $m < \infty$ such that
\begin{align}
\min_{z \in \rmB_r} \min_{N \geq 1} \sfP_z^N \big( \sfH_N(0) \leq m \big) \geq 1/2 \,.
\end{align}
Conditioning on $\sfX_N \big(\sfH_N(\rmB_r(y)) \big)$ and using the strong Markov property, we have
\begin{align}
\sfP_x^N \big( \sfH_N(y) < \vart_N(2n_0) \big) \geq
\sfP_x^N \big( \sfH_N(y) < \vart_N(n_0) + m \big) \geq 
\tfrac12 \sfP_x^N \big( \sfH_N( \rmB_{r} (y) < \vart_N(n_0) \big)  \,,
\end{align}
for all $N$ large enough, depending on $r$ and $n_0$, and any $x,y$ as in the supremum in~\eqref{eq:A.1_jump_control_true}. But~(6.18) from~\cite{JLT} (which corresponds to the case $r = 1$) says that the left hand side above will be smaller than $\epsilon/2$ for all such $x$ and $y$, provided we choose $n_0 > 0$ small enough and $N$ large enough, both depending only on $\epsilon$. The result follows.
\end{proof}

We can now prove:
\begin{lem} Let $\eps >0$, $M \in\bbN$ and $n \in \bbN$. There is $\upsilon(\e,M,n) >0$, $r_0(\eps,M) >0$ and $N_0(\eps,M,r) \in \bbN$ so that when $r \geq r_0$ and $N \geq N_0$, we have
\begin{align}
\sfJ_N^{(r,M)} (\vart_N(n)) \leq \upsilon
\end{align}
with probability at least $1 -\eps$. 
\label{lem:jump_control}
\end{lem}

\begin{proof} For $\eps >0$ and $M \in \bbN$ given, use Proposition \ref{prop:extreme_sep} to choose $r_0 >0$ large so that when $r \geq r_0$ and $N$ is sufficiently large, we have $\bbP( \mathcal{S}_N(r,M) ) \geq 1-\eps/2$. Then, let $\delta >0$ and use Lemma~\ref{lem:jump_control_true} to find $n_0 >0$ and $N_0 > 0$, so that whenever $N \geq N_0$ the statement in \eqref{eq:A.1_jump_control_true} holds with $\delta/M$ in place of $\epsilon$. Using the strong Markov property and the union bound we shall have
\begin{align}
\bbP \left(  \sfR_{k+1} - \sfS_k \geq \vart_N(n_0) \big| \mathcal{S}_N(r,M) \right) \geq 1 - \delta \,,
\label{eq:jump_control_1}
\end{align}
for all $k\geq1$. Now, for $m \in \bbN$, define the event $\mathcal{E}_k$ as
\begin{align}
\mathcal{E}_k := \big\{ \exists j \in \{0,\dots, m-1\} \text{ such that } \sfR_{mk + j +1} - \sfS_{mk +j} \geq \vart_N(n_0) \big\} \,,
\end{align}
and note that iterating the Markov property and using \eqref{eq:jump_control_1} gives, for each $k \geq 1$, the bound $\bbP( \mathcal{E}_k^\rmc | \mathcal{S}_N(r,M) ) \leq (M\delta)^m$. Thus,
\begin{align}
\bbP\left( \mathcal{E}_1 \cap \dots \cap \mathcal{E}_{ \lceil n / n_0 \rceil } \big| \mathcal{S}_N(r,M) \right) \geq 1- \lceil n / n_0 \rceil (M\delta)^m \,.
\end{align}
Observing that $\mathcal{E}_1 \cap \dots \cap \mathcal{E}_{ \lceil n / n_0 \rceil } \subset \big\{ \sfJ_N^{(r,M)} ( \vart_N(n)) \leq 2m \lceil n / n_0 \rceil \big\}$, it follows that
\begin{align}
\bbP \left( \sfJ_N^{(r,M)} (\vart_N(n)) \leq 2m\lceil n / n_0 \rceil \right) \geq \bbP \left(\mathcal{E}_1 \cap \dots \cap \mathcal{E}_{ \lceil n / n_0 \rceil } \big| \mathcal{S}_N(r,M)) \bbP( \mathcal{S}_N(r,M) \right) \,.
\end{align}
It remains to choose $m$ sufficiently large so that the right hand side above is at least $1-\eps$.  \end{proof}

\subsection{Closeness of the random walk to its trace}
\label{ss:TraceClose}

We are finally ready to prove that the trace process and the random walk process are close. Our strategy is similar to the one employed in Section~\ref{sec:CloseK}. Setting $\Delta^{(r,M)}(n) := t(n) - t^{(r,M)}(n)$ to be the size of the ``time lag" between $X_N$ and $X_N^{(r,M)}$ after $n$ steps of $\sfX_N$, the first lemma shows that for $M$ large enough, $\Delta^{(r,M)}(n)$ is small in the scale considered.

\begin{lem} Let $\beta > \alpha$, let $\eps >0$ and let $n \in \bbN$. There is $M_0(\beta,\eps,n) \in \bbN$, $r_0(\beta, \eps,M,n) >0$ and $N_0(\beta, \eps, M,n,r)$ so that $M \geq M_0$, $r \geq r_0$ and $N \geq N_0$ implies
\begin{align}
s_N^{-1}\Delta^{(r,M)} (\vart_N(n)) \leq \eps \,
\end{align}
occurs with probability at least $1- \eps$. 
\label{lem:2_delta_control}
\end{lem}

\begin{proof}
Write $\vart_N(n)$ as $\vart$ throughout the proof, and recall that $\ol{\Lambda}^\rmc_N(r,M) = V_N^* \setminus  \ol{\Lambda}_N(r,M)$. For any realization of the field $h_N$, it follows from  \eqref{eq:discrete_trace_clock_process} that
\begin{align}
\Delta^{(r,M)} ( \vart) = \sum_{j=0}^{\vart-1} \exp\big(\beta h_{N, \sfX_N(j)}\big) \1_{\left\{ \sfX_N(j) \in \ol{\Lambda}_N^\rmc(r,M)\right\}} E_{N,j} \,,
\end{align}
where we recall the $E_{N,j}$ are i.i.d. mean one exponentially distributed random variables independent of $h_N$ and $\sfX_N$. Conditional on $h_N$, we then have: 
\begin{align}
\bbE \left( \Delta^{(r,M)} ( \vart ) \Big| h_N \right) =\sum_{ x \in \ol{\Lambda}_N^\rmc(r,M)} \exp(\beta h_{N,x}) \sfG_{\vart}^N (0,x)\,,
\end{align}
where we recall that $\sfG_n^N(0,x)$ is the expected number of visits to $x$ in $n$ steps for the random walk $\sfX_N$ stating from $0$.

We take $N$ sufficiently large so that the estimate of Lemma \ref{lem:2_green_estimate_one} is valid to obtain:
\begin{align}
s_N^{-1}\bbE\left( \Delta^{(r,M)} ( \vart ) \Big| h_N \right) \leq g^{-1}Cn \sum_{ x \in \ol{\Lambda}_N^\rmc(r,M)} \exp(\beta (h_{N,x} - m_N)) \,,
\label{eq:delta_control_1}
\end{align}
For $\eps >0$, define the event $\mathcal{E}$ which depends only on the field $h_N$:
\begin{align}
\mathcal{E} := \left\{ \sum_{x \in \ol{\Lambda}_N^\rmc(r,M)} \exp(\beta (h_{N,x} -m_N) ) < \frac{g \eps^2}{2Cn} \right\}\,,
\label{eq:4.30}
\end{align}
where $C$ in \eqref{eq:4.30} is as in \eqref{eq:delta_control_1}. Using Proposition~\ref{prop:tightness}, take $M$, then $r$ and then $N$ large enough depending on $n$, $\eps$ and $\beta$ so that $\bbP( \mathcal{E}) \geq 1-\eps/2$. We find
\begin{align}
\bbP\Big( s_N^{-1}\Delta^{(r,M)} (\vart) > \eps \Big) \leq \eps/2 + \eps^{-1} \bbE \big( \bbE [ s_N^{-1} \Delta^{(r,M)} (\vart) | h_N ] \1_{\mathcal{E}} \big)
\leq \eps \,,
\end{align}
where we have used \eqref{eq:delta_control_1}. \end{proof}

In what follows, abbreviate $\wt{t}(n)=t^{(r,M)}(n)$. The next lemma shows that $[0,\infty)$ is partitioned by the intervals $\big[ \wt{t} (\sfR_k), \wt{t}(\sfS_k) \big)$, and that for most of each interval, the processes $X_N$ and $X_N^{(r,M)}$ are not too far apart.

\begin{lem} For each $r >0$ and $M \in \bbN$, the intervals $\{ [ \wt{t}(\sfR_k), \wt{t}(\sfS_k ) ) \}_{k=1}^\infty$ form a partition of $[0,\infty)$. Moreover, for each $t$ in an interval of the form
\begin{align}
\big[ \wt{t}(\sfR_k) , \wt{t}(\sfS_k) \big) \setminus \big[ \wt{t}(\sfR_k) , \wt{t}(\sfR_k) + \Delta^{(r,M)}( \sfR_{k}) \big] \,,
\end{align}
we have $\big\| X_N(t) - X_N^{(r,M)}(t) \big\| \leq 2r_N$.
\label{lem:lag_closeness}
\end{lem}

\begin{proof} By construction, $\wt{t}$ is constant when $\sfX_N(n) \notin \ol{\Lambda}_N(r,M)$, from which it follows that $\wt{t}(\sfS_k) = \wt{t}(\sfR_{k+1})$ for each $k \geq 0$. This settles the first claim. 

Turning to the second claim, observe that $X_N$ and $X_N^{(r,M)}$ are at distance of at most $r_N$ on the intersection $[t(\sfR_k), t(\sfS_k)) \cap [\wt{t}(\sfR_k), \wt{t}(\sfS_k))$. We complete the proof by noting $t(\sfR_k) - \wt{t}(\sfR_k)$ is identically $\Delta^{(r,M)}(\sfR_k)$, and that $t(\sfS_k) \geq \wt{t}(\sfS_k)$.  \end{proof}

The next lemma is the analogue of Lemma~\ref{lem:K_pre_K}; we show that for $r$ and $M$ chosen well and $N$ sufficiently large set of times at which $X_N^{(r,M)}$ and $X_N$ are far has small Lebesgue measure.

\begin{lem} Let $\beta >\alpha, \epsilon >0$ and  $\mathfrak{t} > 0$. There are $M_0(\beta,\eps, \mathfrak{t}) \in \bbN$, $r_0 (\beta, \eps, M, \mathfrak{t}) > 0$ and $N_0(\beta, \eps, M, r, \mathfrak{t}) \in \bbN$ so that whenever $M \geq M_0$, $r > r_0$ and $N \geq N_0$, the event
\begin{align}
\Leb \left( t \in [0,s_N\mathfrak{t}] : \big\| X_N( t) - X_N^{(r,M)} (t) \big\| > 2r_N \right)  \leq \epsilon s_N\,
\label{eq:walk_trace_times_0}
\end{align}
occurs with probability at least $1- \eps$.
\label{lem:walk_trace_times}
\end{lem}

\begin{proof} Let $\eps >0$. As $\mathfrak{t}$ is fixed, use Lemma \ref{lem:2_elapsed_cover} to choose $n = n(\beta, \eps, \mathfrak{t}) \in \bbN$ sufficiently large so that with probability at least $1-\eps$, the event $\mathcal{E} := \{ \wt{t}( \vart_N(n)) > s_N\mathfrak{t}\}$ occurs for all $r >0$ and $M \in \bbN$, provided $N$ is sufficiently large. Treat this $n$ as fixed for the remainder of the proof and write $\vart$ for  $\vart_N(n)$. 

Write $B_N^{(r,M)}(\mathfrak{t})$ for the set in \eqref{eq:walk_trace_times_0}, which plays the same role as the set of times $B_M(\mathfrak{t})$ defined in \eqref{eq:K_pre_K_bad_times}. Appealing to the strategy used in the proof of Lemma \ref{lem:K_pre_K}, let $M, M' \in \bbN$ with $M > M'$, to be determined later and consider the intervals:
\begin{align}
\left\{ I_k \right\} _{k=0}^{\sfJ^{(r,M)}(\vart)+1} := \Big\{ \left[\wt{t}(\sfR_k), \wt{t}(\sfS_k) \right) \Big\}_{k=0}^{\sfJ^{(r,M)} (\vart) +1} \,,
\label{eq:interval_collection}
\end{align}
which by Lemma \ref{lem:lag_closeness} and the definition \eqref{eq:macro_jump_def} of $\sfJ^{(r,M)}(\vart)$ are disjoint and cover the interval $[0, \wt{t}(\vart)]$. Call an interval $I_k$ in this collection \emph{extremely deep} (ED) if $\sfX_N(\sfR_k) \in \ol{\Lambda}_N(r,M')$, and call the interval \emph{moderately deep} (MD) otherwise. By Lemma \ref{lem:lag_closeness}, and within the event $\mathcal{E}$, we find
\begin{align}
B_N^{(r,M)}(\mathfrak{t}) \subset \left( \bigcup_{I_k \text{ is ED} } [ \wt{t}(\sfR_k), \wt{t}(\sfR_k) + \Delta^{(r,M)}(\sfR_k) ] \right) \cup \left( \bigcup_{I_k \text{ is MD}} I_k \right) \,.
\end{align}
Observe that on $\mathcal{S}_N(r,M)$, the number of ED intervals in the collection \eqref{eq:interval_collection} is at most $\sfJ^{(r,M')}(\vart)$ and the total length of all MD intervals in \eqref{eq:interval_collection} is at most $\Delta^{(r,M')}(\vart)$. Thus, on the event $\mathcal{E}\cap \mathcal{S}_N(r,M)$ we have,
\begin{align}
\Leb\big(B_N^{(r,M)}(\mathfrak{t})\big) \leq \sfJ^{(r,M')}(\vart) \Delta^{(r,M)}(\vart) +  \Delta^{(r,M')}(\vart) \,.
\label{eq:walk_trace_times_bound}
\end{align}

We now calibrate parameters. Use Lemma \ref{lem:2_delta_control} to choose $M' = M'(\eps, \beta,n)$ and then $r = r(\eps,\beta,M,n)$ sufficiently large so that $\mathcal{E}_1 := \{ s_N^{-1} \Delta^{(r,M')} (\vart) \leq \eps \}$ occurs with probability at least $1- \eps$ when $N$ is sufficiently large. Next use Lemma \ref{lem:jump_control} and choose $r$ larger if necessary, so that $\mathcal{E}_2 := \{ \sfJ_N^{(r,M')}(\vart) \leq \upsilon\}$ happens with probability at least $1- \eps$ for $N$ sufficiently large.
 Note that $\upsilon$ depends only on the parameters $(\eps, M', n)$, and recall that $n$ has been fixed. Finally, apply Lemma \ref{lem:2_delta_control} once more to choose $M$ even larger than $M'$ so that $\mathcal{E}_3 := \{ s_N^{-1} \Delta^{(r,M)} (\vart) \leq \eps /\upsilon \}$ occurs with probability at least $1- \eps$. We may take $r$ larger if necessary so that $\bbP( \mathcal{S}_N(r,M) ) \geq 1- \eps$. For all $N$ sufficiently large then, by using \eqref{eq:walk_trace_times_bound}, we have 
\begin{align}
\Leb \big(B_N^{(r,M)} (\mathfrak{t}) \big) \leq 2 \eps s_N
\end{align}
on the intersection $\mathcal{E} \cap \mathcal{E}_1 \cap \mathcal{E}_2 \cap \mathcal{E}_3 \cap \mathcal{S}_N(r,M)$, and hence with probability at least $1- 5\eps$. This completes the proof.\end{proof}

We are finally ready for:
\begin{proof}[Proof of Proposition~\ref{prop:TraceClose}]
Thanks to the bounded diameter of $V^*$, the proposition is an immediate consequence of Lemma~\ref{lem:walk_trace_times}.
\end{proof} 

\section{Trap Hopping Dynamics of the Trace Process}\label{s:trap_hop}

As we have seen from the analysis in the previous section, Subsection~\ref{sec:macroscopic_jumps} in particular, the trace process spends most of its time at deep traps, making macroscopic jumps between deep traps almost instantaneously. In this section we show that, scaled properly, the accumulated time in each visit to a deep trap converges to an exponentially distributed random variable whose mean is the depth of the trap. Moreover, at each of the macroscopic jumps, the next trap to visit will be chosen approximately uniformly. We will show that these exponential and uniform random variables are in fact independent, thus the trace process bears the trap hopping characteristics of a spatial pre K-process, as constructed in Section~\ref{sec:CloseK}. The end product of this section, namely Proposition~\ref{prop:close_all}, will be a joint asymptotic description of both this trap hopping behavior and the underlying trapping landscape.  

\subsection{Simple random walk estimates}
We start by providing the needed random walk estimates for $\sfX_N$. Recall the notation 
$\sfP_x^N$, $\sfE_x^N$ for the probability measure and expectation when $\sfX_N(0) = x$ and $\sfH_N(A)$ and $\sfH_N(x)$ for the hitting time of the set $A$ and vertex $x$, as introduced in Subsection~\ref{sec:macroscopic_jumps}
 
The first lemma, which we borrow almost verbatim from~\cite{JLT}, concerns the hitting measure of separated balls in $V_N^*$. This lemma will be used to show that deep traps are selected almost uniformly at each macroscopic jump of $\sfX_N$. Specifically, we consider $A_N \subset V_N^*$ such that 
\begin{enumerate}
\item [(1)] $|A_N| = M$\,,
\item [(2)]$\min_{x \neq y \in A_N} d_N^*(x,y) > r_N /2 $\,,
\end{enumerate}
Then,
\begin{lem} Let $\eps >0$, $M \in \bbN$ and $r >0$. Let $A_N = \{ x_{N,1}, \dots, x_{N,M} \}$ satisfy (1) and (2) directly above. There is $N_0(r,\eps,M) \in \bbN$ so that whenever $N \geq N_0$,
\begin{align}
\label{e:5.1}
\sup_{y \in \big(\rmB_{r_N/2}(A_N)\big)^\rmc} \left| \sfP_y^N \left[ \sfH_N(\rmB_r(x_{N,1})) < \sfH_N(\rmB_r(A_{N,1})) \right]  - \frac{1}{M} \right| < \eps \,,
\end{align}
where we write $A_{N,i}$ for the $A_N \setminus \{x_{N,i} \}$. 
\label{lem:ord_close_pre}
\end{lem}

\begin{proof}
We sketch changes that need to be made to the proof of Lemma 6.9 in \cite{JLT}. Key inputs to this proof are (6.15), (6.16) and (6.18) in \cite{JLT}. To prove Lemma~\ref{lem:ord_close_pre}, we use (6.15) unchanged, but we replace (6.16) and (6.18) with corresponding estimates for hitting times of $\rmB_r(x_{N,i})$ in place of $x_{N,i}$. These estimates follow from (6.16) and (6.18) immediately via a union bound, and it is here that $N_0$ inherits its dependence on $r$. The rest of the proof goes through, again replacing instances of $x_{N,i}$ by $\rmB_r(x_{N,i})$. \end{proof}

The next estimate concerns the scaling limit of the time spent by the walk inside a ball of radius $r$, before exiting the ball of radius $r_N$. This will be used to show that the time spent visiting a deep trap scales in law to an exponential random variable. Here, we lose no generality working with a simple random walk $\sfX$ on $\bbZ^2$, writing law and expectation as $\sfP_x$ and $\sfE_x$.

Let $\sfH^{r_N}$ be the first time $\sfX$ exits the ball $\rmB_{r_N} \subset \bbZ^2$ and for $y \in \rmB_r$, define
\begin{align}
\sfL_N(y) := \# \{ n \in [0, \sfH^{r_N} ] : \sfX(n) = y \}\,,
\end{align}
to be the the number of visitis of $y$ before step $\sfH^{r_N}$. Let us also define
\begin{equation}
L_N(y) := \sum_{k=1}^{\sfL_N(y)} E_{j} \,,
\label{e:5.3}
\end{equation}
where $(E_{j} : j \geq 1)$ are i.i.d. mean one exponential random variables, independent of $\sfX$, exactly as in Subsection~\ref{ss:RWClock}. The random variable $L_N(y)$ can be thought of as the local time of a continuous version of $\sfX$ with unit mean holding times, and we shall refer to $L_N(y)$ as such.

The following lemma shows that the $(L_N(y))_{y\in \rmB_r}$ converge weakly under proper scaling to the same exponential random variable. Below we write $\ol{1}_r$ for a vector of ones indexed by $\rmB_r$.

\begin{lem}  
Let $r > 0$ and $x \in \rmB_r$. Then under $\sfP_x$,
\begin{equation}
\label{e:5.4}
L_N / \log N \equiv \big(L_N(y) / \log N : y \in \rmB_r\big) \underset{N \to \infty}{\Longrightarrow}
 e \ol{1}_r 
\end{equation}
where $e$ is an exponential random variable with mean $g$.
\label{lem:lt_exp_prelim}
\end{lem}

\begin{proof} Suppose first that $x=0$. In this case $\sfL_N(0)$ is geometrically distributed with parameter $p_N := 1 / \sfG_{\rmB_{r_N}}(0,0)$, where for a set $A \subset \bbZ^2$, we recall that $\sfG_A$ denotes the Green function associated with a simple random walk killed upon exit from $A$. Since $\sfG_{\rmB_{r_N}}(0,0) \sim g \log N$ as $N \to \infty$ (see Lemma~\eqref{lem:GAsymp}), it is elementary that $\sfL_N(0)/\log N$ converges weakly to an exponential random variable $e$ with mean $g$.

Next, we use a result by Auer~\cite{Auer}: if $\xi(y,n)$ denotes the number of visits to vertex $y$ within the first $n$ steps of $\sfX$, then
\begin{align}
\lim_{n \to \infty} \sup_{y \in \rmB_r} \left| \frac{ \xi(y,n)}{ \xi(0,n) } - 1 \right| = 0 \,,
\label{eq:auer_content_1}
\end{align} 
almost surely. Observe that as $\sfH^{r_N} \to \infty$ with $N$, we also have
\begin{align}
\label{e:5.6}
\lim_{N \to \infty} \sup_{y \in \rmB_r} \left| \frac{ \sfL_N(y)}{ \sfL_N(0) } - 1 \right| = 0 \,,
\end{align}

By recurrence, $\sfL_N(y) \to \infty$ as $N \to \infty$ almost surely, and the strong law of large numbers implies $L_N(y)/\sfL_N(y) \to 1$ for all $y \in \rmB_r$ almost surely. It follows that~\eqref{e:5.6} holds with $\sfL_N(y)$ replaced by $L_N(y)$. Dividing then both the numerator and denominator in~\eqref{e:5.6} by $\log N$ and using the weak convergence of $\sfL_N(0)/\log N$ to $e$, we obtain~\eqref{e:5.4} when $x=0$.

For general $x \in \rmB_r$, let $\sfH(0)$ be the hitting time of $0$, and write $\sfL_N(y)$ as the sum $\sfL^{(1)}_N(y) + \sfL^{(2)}_N(y)$, where $\sfL^{(1)}_N(y)$ is the number of visits to $y$ before step $\sfH^{r_N} \wedge \sfH(0)$, while $\sfL^{(2)}_N(y)$ is the remaining number of visits to $y$ until time $\sfH^{r_N}$. Let $L^{(1)}_N(y)$ and $L^{(2)}_N(y)$ be the corresponding local times; recurrence of $\sfX$ and the Markov property show that
\begin{equation}
\sup_{y \in \rmB_r} L^{(1)}_N(y)/(\log N)
\leq \frac{1}{\log N} \sum_{k=1}^{\sfH^{r_N} \wedge \sfH(0)} E_k \underset{N \to \infty} \longrightarrow 0 
\text{ a.s.}\quad\,,\quad L^{(2)}_N/\log N \underset{N \to \infty}{\Longrightarrow} e \ol{1}_r\,.
\end{equation}
This shows~\eqref{e:5.4} in the general case and completes the proof.
\end{proof}

\subsection{Trap hopping dynamics at large $N$}

We now use the estimates from the previous subsection to describe the limiting joint law of the indices of visited traps and the local time spent near these traps. The end product of this subsection, namely Proposition~\ref{prop:close_all}, will be the convergence in law for these variables together with the trapping landscape of the field.

Recall that $\cS_N(r,M)$ from~\eqref{eq:well-separated} denotes the event that $h_N$ is $(r,M)$-separated and that the stopping times $\sfR_k, \sfS_k$ from Subsection~\ref{sec:macroscopic_jumps} mark the beginning and end of visits to a deep trap. For each $k\geq1$, when $\cS_N(r,M)$ occurs, $\sfX_N(\sfR_k)$ is closest to a unique trap in $\ol{\Lambda}_N(r,M)$ and we let $I_{N,k}'$ be the index of this trap. Formally, $I_{N,k}'$ is the unique {\em ordinal} in $\{1, \dots, M\}$ such that $\sfX_N(\sfR_k) \in \rmB_r(x_{N,I_{N,k}'})$, where we recall that
$x_{N,1}, \dots, x_{N,M}$ are the deep traps listed in descending order of their depth, per~\eqref{eq:Lambda_def}.

If $h_N$ is not $(r,M)$-separated, $I_{N,k}'$ might not be well-defined, and for this purpose we take $(\wt{U}_{N,k})_{k \geq 1}$ to be a collection of independent uniform random variables on $\{1, \dots, M\}$ which are also independent of everything else, and set
\begin{align} 
I_{N,k} := I_{N,k}' \1_{ \mathcal{S}_N(r,M)} + \wt{U}_{N,k} \1_{\mathcal{S}_N(r,M)^\rmc} \,.
\label{eq:true_ordinal_def}
\end{align}
As for the local time spent during the $k$-th visit to a trap, we again assume first that $\cS_N(r,M)$ occurs and, in this case, for all $y \in \rmB_r$ define:
\begin{align}
L_{N,k}'(y) := \sum_{j =\sfR_{k}}^{\sfS_{k}} \1_{\{\sfX_N(n) = x_{N,I_{N,k}} + y\}} E_{j} \,,
\end{align}
where the $E_{j}$ are as in the previous subsection. We then take in addition a collection $(\wt{e}_k)_{k \geq 1}$ of i.i.d. mean $g$ exponential random variables, independent of everything else and set:
\begin{align}
\ell_{N,k}(y) := \big(L_{N,k}'(y)/\log N\big) \1_{\mathcal{S}_N(r,M)} + \wt{e}_k 
\ol{1}_r \1_{\mathcal{S}_N(r,M)^\rmc} \,.
\label{eq:true_loc_def}
\end{align}
We shall write $\ell_{N,k}$ for the collection $(\ell_{N,k}(y) :\: y \in \rmB_r)$.

The following lemma shows that the joint law of $(I_{N,k})_{k \geq 1}$ and $(\ell_{N,k})_{k \geq 1}$ for large $N$ is approximately that of independent uniform and exponential random variables. 
\begin{lem}
\label{lem:independence_day2}
Let $(e_k)_{k \geq 1}$ be i.i.d. exponential random variables with mean $g$ and let $(U_k)_{k\geq 1}$ be i.i.d. uniform random variables on the set $\{1, \dots, M\}$, with both collections independent of each other. For $K \in \bbN$, let $A = \prod_{k=1}^K \big(A_k \times \{u_k\}\big)$ where $u_1, \dots, u_K \in \{1, \dots, M\}$ and $A_1, \dots, A_K \subseteq \bbR^{\rmB_r}$ are $A_k$ measurable and continuous with respect to the law of $e_k \ol{1}_r$. Then, 
\begin{equation}
\label{eq:144}
\bbP \Big( \big(\ell_{N,k} ,\, I_{N,k} \big)_{k=1}^K \in A \,\Big|\, h_N \Big) 
\longrightarrow \bbP \big(\big(e_k \ol{1}_r ,\, U_k \big)_{k=1}^K \in A \big) \,,
\end{equation}
as $N \to \infty$, where convergence takes place in the $\rmL^\infty$-norm on the underlying probability space.
\end{lem}

\begin{proof} 
Treat the parameters $(r,M)$ as fixed in this proof. By the definitions of $I_{N,k}$ and $\ell_{N,k}$, on the event $\cS_N(r, M)^\rmc$ that $h_N$ is not $(r,M)$-separated, both probabilities in~\eqref{eq:144} are equal to each other and hence it suffices to show $\rmL^\infty$-convergence on $\cS_N(r,M)$. For what follows, we shall write $\cF_\sfT$ for the sigma-algebra generated by the stopping time $\sfT$ with respect to the natural filtration of $\sfX_N$. 

The proof goes by induction on $K$. Recall that $\sfR_k$ and $\sfS_k$ are the entrance and exit times of the $k$-th trap visited by $\sfX_N$. For $K=1$, condition on $\cF_{\sfR_1}$ and use the strong Markov property to write the probability $\bbP \big(\ell_{N,1} \in A_1,\, I_{N,1} = u_1 \,\big| h_N \big)$ as
\begin{equation}
\label{eq:141}
\bbE \Big( \1_{\{I_{N,1} = u_1\}} \sfP_{\sfX_N(\sfR_1) - x_{N,u_1}} \big(L_N/\log N \in A_1\big) \, \Big| h_N \Big) \,,
\end{equation}
where $L_N$ is as in \eqref{e:5.3}, and where $\sfX_N(\sfR_1) - x_{N,u_1}$ is thought of as a vertex in $\rmB_r \subset \bbZ^2$. Writing $\Lambda_{N, u_k}$ for the set $\Lambda_N(r,M) \setminus \{x_{N,u_k}\}$, we have
\begin{align}
\bbP \big(I_{N, 1} = u_1 \,\big|\, h_N\big) = \sfP^N_0 \big( \sfH_N(\rmB_r(x_{N,u_1})) < \sfH_N(\rmB_r(\Lambda_{N,{u_1}}))\big)\,.
\label{eq:142}
\end{align}
Since all traps are at least a distance of $r_N/2$ from $\sfX_N(0) \equiv 0$ when $\cS_N(r, M)$ occurs and $N$ is large enough, Lemma~\ref{lem:ord_close_pre} shows that the above probability converges to $\bbP(U_1 = u_1) = 1/M$ as $N \to \infty$ under the $\rmL^\infty$-norm. At the same time, by Lemma~\ref{lem:lt_exp_prelim} and the fact that 
$\sfX_N(\sfR_1) - x_{N,u_1} \in \rmB_r$ on $\{I_{N,1} = u_1\}$, which is a finite set, the probability in~\eqref{eq:141} converges to 
$\bbP(e_1 \ol{1}_r \in A_1)$ as $N \to \infty$ in the $\rmL^\infty$-sense on $\{I_{N,1} = u_1\}$. It follows from H\"{o}lder's inequality under the conditional measure $\bbP(\cdot|h_N)$ that~\eqref{eq:141} converges to $\bbP(U_1 = u_1)  \bbP(e_1 \ol{1}_r \in A_1)$
in $\rmL^\infty$ on $\cS_N(r,M)$. This settles the base case. 

Let $K \geq 2$ and assume the statement of the lemma holds with $K-1$ in place of $K$. Write $A = A' \times \big(A_{K} \times \{u_{K}\}\big)$ with $A' = \prod_{k=1}^{K-1} \big(A_k \times \{u_k\}\big)$, and let $\cA'_N$ be the event
$\big\{\big(\ell_{N,k} ,\, I_{N,k} \big)_{k=1}^{K-1} \in A'\big\}$. Conditioning on $\cF_{\sfS_{K-1}}$ and $\cF_{\sfR_K}$ and using the strong Markov property, the left hand side of~\eqref{eq:144} is equal to
\begin{equation}
\label{eq:145}
\bbE \Big( \1_{\cA'_N} \sfE^N_{\sfX_{\sfS_{K-1}}} \Big( 
\1_{\{\sfH_N(\rmB_r(x_{N,u_K})) < \sfH_N(\rmB_r(\Lambda_{N,{u_K}}))\}}
\sfP_{\sfX_N(\sfR_K) - x_{N,u_K}} \big(L_N/\log N \in A_K\big) \Big) \, \Big| h_N \Big).
\end{equation}

Since $\sfX_{\sfS_{K-1}}$ is by definition at least $r_N/2$ away from all traps, it follows as before, that the middle expectation in~\eqref{eq:145} converges in $\rmL^\infty$-sense to $\bbP(U_K=u_K) \bbP(e_K \ol{1}_r \in A_K)$ as $N \to \infty$ on $\cS_N(r,M)$. On the other hand, the induction hypothesis gives that $\bbP(\cA'_N | h_N)$ converges to 
$\bbP \big((e_k \ol{1}_r ,\, U_k)_{k=1}^{K-1} \in A' \big)$ as $N \to \infty$ again in $\rmL^\infty$. Using these together with H\"{o}lder's inequality for $\bbP(\cdot|h_N)$ shows that~\eqref{eq:145} converges as $N \to \infty$ to the product
\begin{align}
\bbP \left(\big(e_k \ol{1}_r ,\, U_k \big)_{k=1}^{K-1} \in A' \right) 
\bbP (U_K = u_K) \bbP(e_K \ol{1}_r \in A_K)
\end{align}
under the $\rmL^\infty$-norm. This is precisely the limit in~\eqref{eq:144}.
\end{proof}

Finally, we treat the asymptotic joint law of both the trap hopping dynamics of the trace process and the trapping landscape of the underlying field.
\begin{prop}
\label{prop:close_all}
Fix $\beta > \alpha$ and let $M \in \bbN$, $K \in \bbN$. The joint law of 
\begin{align}
\label{eq:7.2.1}
\big(\ell_{N,k} ,\, I_{N,k} \big)_{k=1}^K
\ , \quad\quad
\big( x_{N,j}/N ,\, \rme^{-\beta m_N} \tau_r(x_{N,j}) \big)_{j=1}^M 
\end{align}
converges weakly as $N \to \infty$ followed by $r \to \infty$ to the joint law of 
\begin{align}
\label{eq:7.2.2}
\big(e_k \ol{1}_\infty ,\, U_k \big)_{k=1}^K
\ , \quad\quad
\big( \xi_j ,\, \tau_j \big)_{j=1}^M \,,
\end{align}
where $\big((e_k, U_k))_{k=1}^K$ are as in Lemma~\ref{lem:independence_day2} and $\big((\xi_j, \tau_j)\big)_{j=1}^M$ are as in Proposition~\ref{prop:depth_close}, with both collections independent of each other, and where $\ol{1}_\infty$ is a vector of ones indexed by $\bbZ^2$.
\end{prop}

\begin{proof}
Let $r > 0$ and take $A$ and $B$ to be any measurable continuity sets for $\big(e_k \ol{1}_r ,\, U_k \big)_{k=1}^K$ and $\big( \xi_j ,\, \tau_j \big)_{j=1}^M$ respectively, with $A$ having the form in Lemma~\ref{lem:independence_day2}. Conditioning on $h_N$ write,
\begin{multline}
\label{eq:149}
\bbP \Big( \big(\ell_{N,k} ,\, I_{N,k} \big)_{k=1}^K \in A ,\,
\big( x_{N,j}/N ,\, \rme^{-\beta m_N} \tau_r(x_{N,j}) \big)_{j=1}^M  \in B \Big) \\ = 
\bbE \Big[ \bbP \Big(\big(\ell_{N,k} ,\, I_{N,k} \big)_{k=1}^K \in A \, \Big|\, h_N\Big)
\1_B\Big(\big(x_{N,j}/N ,\, \rme^{-\beta m_N} \tau_r(x_{N,j}) \big)_{j=1}^M\Big) \Big]
\end{multline}
By Lemma~\ref{lem:independence_day2}, the first term in the expectation goes to
$\bbP \big(\big(e_k \ol{1}_r ,\, U_k \big)_{k=1}^K \in A \big)$ as $N \to \infty$ in the $\rmL^\infty$-norm, for any $r > 0$. 
At the same time, by Proposition~\ref{prop:depth_close} the expectation of the second term 
goes to $\bbP \big(\big( \xi_j ,\, \tau_j \big)_{j=1}^M \in B \big)$ as $N \to \infty$ followed by $r \to \infty$. It follows by H\"older's inequality that the right hand side of~\eqref{eq:149} converges in the stated limits to the product of the last two probabilities. This product is precisely,
\begin{equation}
\bbP \Big(\big(e_k \ol{1}_r ,\, U_k \big)_{k=1}^K \in A
,\, \big( \xi_j ,\, \tau_j \big)_{j=1}^M \in B \Big)  \,.
\label{eq:5.19}
\end{equation}
The proof is completed by observing that the collection of events of the form appearing in \eqref{eq:5.19} is a convergence determining class for the distributions in question.
\end{proof}

\section{Conclusion of the Proof}\label{s:final}
In this section we complete the proof of Theorem~\ref{thm:A}. We first use the results from the previous section to prove that the trace process can be coupled together with the $\chi$-driven spatial pre K-process, so that with high probability, they are close in the $\|\cdot\|_{\rmL([0,\mft], V^*)}$-metric. We then use the closeness of random walk to its trace (Section~\ref{s:TraceClose}) and the closeness of the K-process to the pre K-process (Section~\ref{sec:CloseK}) to complete the proof.

\subsection{Closeness of the trace process and the pre K-process}
Recall that $Y^{(\beta)}_M$ denotes the $\chi$-driven pre K-process, as introduced in Section~\ref{sec:CloseK}. The goal in this subsection is to show that $X^{(r, M)}_N$ and $Y^{(\beta)}_M$ can be coupled so that they are close in the $\|\cdot\|_{\rmL([0,\mft], V^*)}$-metric with high probability, provided we choose $N$ and $r$ appropriately:
\begin{prop}
Fix $\beta > \alpha$ and let $\eps >0$, $M \in \bbN$, and $\mathfrak{t} > 0$. There is $r_0(\eps, M, \mathfrak{t}) > 0$ and $N_0(\eps,M, r, \mathfrak{t}) \in \bbN$ so that when $r \geq r_0$ and $N \geq N_0$ there is a coupling between $X_N$ and $Y_M^{(\beta)}$ so that
\begin{align}
\Big\| \frac{1}{N} X^{(r, M)}_N(s_N\cdot) - Y^{(\beta)}_M(\cdot) \Big\|_{\rmL([0, \mathfrak{t}], V^*)} < \epsilon \,
\end{align}
hold with probability at least $1-\eps$. 
\label{prop:trace_pre_K}
\end{prop}

Let $\cC_M$ denote the limiting objects in Proposition~\ref{prop:close_all}, namely the random variables in~\eqref{eq:7.2.2}. To prove Proposition~\ref{prop:trace_pre_K}, we first show that there is a coupling between $(X_N, h_N)$ and $\cC_M$ so that the time spent by the trace process during visit $k$ to trap $j$ is close to $e_k \tau_j$.

\begin{lem} Let $\eps >0$, $K \in \bbN$ and $M \in \bbN$. There is $r_0(\eps,K,M)$ and $N_0(\eps, K,M, r)$ so that whenever $r \geq r_0$ and $N \geq N_0$, there is a coupling of $(X_N, h_N)$ with $\mathcal{C}_M$, so that both
\begin{align}
\sum_{k=1}^K \left|s_N^{-1} \big| t^{(r,M)} (\sfS_k) - t^{(r,M)}(\sfR_k) \big| - e_k \tau_{U_k} \right| \leq \eps 
\label{eq:first_coupling_statement}
\end{align}
and
\begin{equation}
\label{eq:second_coupling_statement}
\sup_{k \leq K} \sup_{n \in [\sfR_k, \sfS_k)} \big \|\tfrac{1}{N}\sfX_N(n)- \xi_{U_k} \big \| \leq \epsilon
\end{equation}
occur with probability at least $1- \eps$. 
\label{lem:first_coupling}
\end{lem}

\begin{proof} Treat $M$ and $K$ as fixed throughout the proof. For $\delta >0$ and $\epsilon >0$, apply Proposition~\ref{prop:close_all} to find $r_0(\eps,M) >0$ and $N_0(\eps, M, r) \in \bbN$ so that when $r \geq r_0$ and $N \geq N_0$,
we may couple $(X_N, h_N)$ with  $\cC_M$, such that (1) -- (4) below hold with probability at least $1-\epsilon/2$:
\begin{enumerate} 
\label{eq:first_coupling_loc}
\item [(1)] For $k \in \{1, \dots, K\}$ and $y \in \rmB_r$, $ | \ell_{N,k}(y) - e_k | < \delta $.
\item [(2)] For $k \in \{1, \dots, K\}$, $I_{N,k} = U_{k}$. \label{eq:first_coupling_ordinal}
\item [(3)] For $j \in \{1, \dots, M\}$, $\big\|x_{N,j}/ N - \xi_j \big\| < \delta$. \label{eq:first_coupling_position}
\item [(4)] For $j \in \{1, \dots, M\}$, $| \rme^{-\beta m_N} \tau_r( x_{N,j} ) - \tau_j | < \delta $. \label{eq:first_coupling_depth}
\end{enumerate}
Recall from \eqref{eq:discrete_trace_clock_process} that for $n \in \bbN$, $t^{(r,M)}(n)$ denotes the time accumulated by the trace process in $n$ steps of the embedded discrete time simple random walk. Since, by definition, during steps $[\sfR_k, \sfS_k)$ the walk is visiting trap $I_{N,k}$, the total time spent in this trap is
\begin{align}
\mathfrak{I}^{(r,M)}_k :=  t^{(r,M)} (\sfS_k) - t^{(r,M)}(\sfR_k) \equiv \sum_{n = \sfR_k}^{\sfS_k - 1} \rme^{\beta h_{N, \sfX_N(n)}} E_{N,n} \textbf{1}_{\{ \sfX_N \in \rmB_r( x_{N, I_{N,k} } ) \}}  \,.
\end{align}
Divide through by $s_N$ and decompose the sum defining $\mathfrak{I}^{(r,M)}_k$ using the local times in \eqref{eq:true_loc_def}: 
\begin{align}
s_N^{-1} \mathfrak{I}_k^{(r,M)} = \sum_{y \in \rmB_r} \ell_{N,k}(y)\rme^{-\beta m_N} \rme^{\beta h_{N, y + x_{N, I_{N,k} } } }  \,.
\label{eq:first_coupling_1}
\end{align}
Note that, for equality to hold directly above, we must work within the high probability event $\mathcal{S}_N(r,M)$ from Proposition~\ref{prop:extreme_sep} that $h_N$ is $(r,M)$-separated, which we can guarantee by taking $r$ and then $N$ larger if necessary. 

Thus, using (1),(2) and (4) in \eqref{eq:first_coupling_1}, we find
\begin{align}
\sum_{k=1}^K | s_N^{-1} \mathfrak{I}_k^{(r,M)} - e_k \tau_{U_k} | \leq K\delta( \delta + \tau_1) + \delta \sum_{k=1}^Ke_k \,,
\end{align}
where we have used that $\tau_1 \geq \tau_j$ for all $j \geq 1$. 
Finiteness of $\tau_1$ and $e_1, \dots, e_K$ almost surely imply that we may further find $L=L(\epsilon, K) > 0$ such that with probability at least $1- \eps$,
\begin{align}
\sum_{k=1}^K | s_N^{-1} \mathfrak{I}_k^{(r,M)} - e_k \tau_{U_k} | \leq L\delta  \,.
\end{align}
Letting $\delta = \epsilon/L$ implies~\eqref{eq:first_coupling_statement}. At the same time, it follows by definition that on $[\sfR_k, \sfS_k)$ we have $\big \|\sfX_N(n) - x_{N, I_{N,k}} \big\| \leq r_N$. Thus on the event that (2) and (3) hold, we also have~\eqref{eq:second_coupling_statement} for all $N$ large enough.
\end{proof}

We are now ready for:
\begin{proof}[Proof of Proposition~\ref{prop:trace_pre_K}]
Fix $M$, $\beta$, $\epsilon$ and $\mft$. Begin by using Lemma \ref{lem:2_elapsed_cover} and Lemma \ref{lem:jump_control} to find $K \in \bbN$ so that with probability at least $1-\eps/2$, we have $s_N^{-1} t^{(r,M)} (\sfR_{K+1}) > \mathfrak{t}$ for all $r$ and then $N$ large enough. With this $K$ and our fixed $M$, increasing $r$ and $N$ if necessary, invoke Lemma~\ref{lem:first_coupling} to find a coupling of $(X_N, h_N)$ and $\cC_M$, under which, with probability at least $1-\epsilon$, both~\eqref{eq:first_coupling_statement} and~\eqref{eq:second_coupling_statement} hold with $\epsilon$ replaced by $\epsilon/(2K(\mft+1))$.

The random variables in $\cC_M$ may be used to build a $\chi$-driven pre K-process. Indeed, take $A = (A(u) : u \geq 0)$ to be a Poisson process on $\bbR_+$ independent of $\mathcal{C}_M$ and having intensity measure $M \rmd u$. For $k \in\{1,\dots, M\}$ construct $A_k = (A_k(u) : u \geq 0)$ by setting 
\begin{align}
A_k(u) = \sum_{j=1}^{A(u)} \1_k(U_j)
\end{align}
Standard Poisson thinning shows $(A_k)_{k=1}^M$ form independent Poisson processes with intensity $1 \rmd u$. Setting also $e_j^{(k)} := e_\ell$ where $\ell$ is the smallest index such that $\sum_{i=1}^\ell \1_{k}(U_i) = j$, it follows that $(e_j^{(k)} : j,k \geq 1)$ are
are i.i.d. mean one exponential random variables, which are also independent of $(A_k)_{k=1}^M$. 
Using $(A_k)_{k=1}^M$ together with $(\xi_j, \tau_j)_{j=1}^M$ from $\cM$ we can define $T_M(\cdot)$, $K_M^{(\tau)}(\cdot)$ and $Y_M^{(\xi, \tau)}(\cdot)$  as in~\eqref{eq:2_pre_clock},~\eqref{eq:pre_K_def} and~\eqref{eq:Ydef} and use these to construct $Y_M^{(\beta)}$ as before. 

Let $B(\mathfrak{t})$ denote the following set of \emph{bad} times:
\begin{align}
B(\mathfrak{t}) := \left\{ t \in [0,\mathfrak{t}] :  \big\| \tfrac{1}{N} X_N^{(r,M)} (s_Nt )- Y_M^{(\beta)}(t) \big\| > \eps/(2 \mft) \right\}\,.
\end{align}
It is immediate from the construction and the choice of coupling that $\Leb(B(\mft)) \leq \epsilon/2$ with probability at least $1-\epsilon$. Recalling that $V^*$ has finite diameter, the proof is complete.\end{proof}

\subsection{Proof of Theorem~\ref{thm:A}}
Fix $\beta > \alpha$, $\mathfrak{t} >0$ and $\eps > 0$.
Using Proposition \ref{prp:K_pre_K}, we find $M_0(\eps, \mathfrak{t}) \in \bbN$ so that $M \geq M_0$ implies that 
\begin{align}
\Big\| Y_M^{(\beta)}(\cdot) - Y^{(\beta)}(\cdot) \Big\|_{\rmL([0,\mathfrak{t}] , V^* ) } < \eps \,
\end{align}
occurs with probability at least $1-\eps$. Take $M$ larger still if necessary, and use Proposition \ref{prop:TraceClose} to deduce that, for $r$ and then $N$ taken sufficiently large, with probability at least $1-\epsilon$,
\begin{align}
\left\| \tfrac{1}{N} X_N (s_N \cdot ) - \tfrac{1}{N} X_N^{(r,M)} (s_N \cdot ) \right\|_{\rmL([0,\mathfrak{t}] , V^* ) } < \eps \,.
\end{align}
Finally, increasing $r$ even more if needed and taking $N$ large enough, we use Proposition~\ref{prop:trace_pre_K} to find a coupling between $(X_N, h_N)$ and $Y_M^{(\beta)}$ under which
\begin{align}
\Big\| \tfrac{1}{N} X^{(r, M)}_N(s_N\cdot) - Y^{(\beta)}_M(\cdot) \Big\|_{\rmL([0, \mathfrak{t}], V^*)} < \epsilon \,
\end{align}
occurs with probability at least $1-\epsilon$.

Altogether, there is a coupling between $(X_N, h_N)$ and $Y^{(\beta)}$ so that for all $N$ large enough, 
\begin{equation}
\Big\| \tfrac{1}{N} X_N(s_N \, \cdot) - Y^{(\beta)}(\cdot) \Big\|_{\rmL([0,\mft], V^*)} < 3\epsilon \,.
\end{equation}
with probability at least $1-3 \epsilon$. Since $\epsilon$ was arbitrary, this shows~\eqref{e:1.8} and completes the proof.

\appendix
\section{Discrete Potential Theory in Two Dimensions}
\label{s:app}
In this section we include bounds for the Green function associated with the simple random walks on the torus $V_N^*$ and on $\bbZ^2$. We denote the latter by $\sfX_N$ and $\sfX$ respectively, and recall that $\sfP_x^N$, $\sfE_x^N$ and $\sfP_x$, $\sfE_x$ respectively denote the underlying probability measure and expectation, when the starting point is $x$. The hitting times of a set $A$ and vertex $y$ are denoted by $\sfH_N(A)$, $\sfH(A)$ and $\sfH_N(y)$, $\sfH(y)$ respectively. Finally, the associated Green functions are defined below:
\begin{equation}
\sfG_n^N(x,y) = \sfE_x^N \sum_{k=0}^n \1_y(\sfX_N(k)) 
\quad , \qquad
\sfG_A(x,y) = \sfE_x \sum_{k=0}^{\sfH(A^{\rmc})-1} \1_y(\sfX(k)) \,,
\end{equation}
for $x,y$ in $V_N^*$ or $\bbZ^2$, $n \geq 0$ and $A \subseteq \bbZ^2$.

We start with standard bounds on $\sfG_{V_N}$, where we recall that $V_N = [0,N)^2 \cap \bbZ^2$. They can be found, for instance, in~\cite{daviaud2006extremes}.
\begin{lem}
\label{lem:GAsymp}
There is an absolute constant $C > 0$ such that for all $N \geq 1$, $x,y \in V_N$,
\begin{equation}
\label{e:GAsymp}
\sfG_{V_N} (x,y) \leq g \log \left( \frac{N}{ \| x -y\| \vee 1} \right)+ N
\end{equation}
Moreover for any $\epsilon > 0$ there is $C=C(\epsilon) > 0$ such that
\begin{equation}
\label{e:GAsymp_LB}
\sfG_{V_N} (x,y) \geq g \log \left( \frac{N}{ \| x -y\| \vee 1} \right) - C \,,
\end{equation}
for all $N \geq 1$ and $x,y \in V_N$ with $\min_{z \in V_N^{\rmc}}\big \{\|x-z\| \wedge |y-z\|\big\} > \epsilon N$.
\end{lem}

Next, we provide bounds on $\sfG_{\vart_N(n)}^N$, recalling that $\vart_N(n) = n \lceil N^2 \log N \rceil$. Although these bounds are standard, we could not find a reference for them and hence their short proof is provided. 
\begin{lem}  There are $C, c \in (0, \infty)$ so that when $N$ is sufficiently large, for all $x \in V_N^*$ and $n \geq 1$ we have
\begin{align}
cn \leq  \frac{\sfG_{\vart_N(n)}^N (0,x)}{\log N} \leq Cn \,.
\end{align}
\label{lem:2_green_estimate_one}
\end{lem}

To prove this lemma, we need an upper bound on $\sfG_{n}^N(0,0)$ for $n$ on the order of $N^2$. 
\begin{lem} Let $\gamma \in \bbN$. There is $c(\gamma) >0$ and $N_0 \in \bbN$ so that $N \geq N_0$ implies
\begin{align}
\sfG_{\gamma N^2}^N (0,0) \leq c \log N\,.
\end{align}
\label{lem:step_green_estimate}
\end{lem}
\vspace{-7mm}

\begin{proof} Let $B_{LN} = [-LN,LN)^2 \cap \bbZ^2$. Consider the simple random walk $\sfX$ on $\bbZ^2$, and first choose $L \in \bbN$ large enough so that the following estimate holds:
\begin{align}
p := \sfP_0( \sfH( \pa B_{LN}) < \gamma N^2 ) < 1/2 \,,
\label{eq:step_green_estimate_1}
\end{align}
where the choice of $L$ depends on $\gamma$, and is fixed henceforth, and where such a choice is possible by Proposition 2.1.2 of \cite{Lawler}, for instance. Through the natural coupling of $\sfX$ and $\sfX_N$, and via the fact that the Green function is maximized on the diagonal, the following bound holds:
\begin{align}
\sfG_{\gamma N^2}^N(0,0) \leq c_L (1-p) \sum_{k \geq 1} k p^k \sfG_{B_{LN}}(0,0)
\end{align}
where the constant $c_L$ is positive, depending only on $L$, and results from the finite number of points mapped to $0 \in V_N^*$ from $B_{LN} \subset \bbZ^2$ by the natural quotient map $\varphi_N  : \bbZ^2 \to V_N^*$.  

Via \eqref{eq:step_green_estimate_1}, 
\begin{align}
\sfG_{\gamma N^2}^N(0,0) &\leq 4c_L  \sfG_{ B_{LN} } (0,0) \\
&\leq 4c_Lg \log L N + 4c_L \wt{c}\,, 
\label{e:A.7}
\end{align}
where \eqref{e:A.7} follows from Lemma~\ref{lem:GAsymp}. As $L$ depends only on $\gamma$, the proof is complete. \end{proof}

\begin{proof}[Proof of Lemma~\ref{lem:2_green_estimate_one}]
Let $\wt{\sfX}_N$ be the lazy simple random walk on $V_N^*$ which waits at each vertex with probability $1/2$. Write $\wt{\sfP}_x$ for the law of $\wt{\sfX}_N$ started at $x \in V_N^*$, suppressing the $N$ in the law for notational convenience. Recall that the \emph{$\epsilon$-uniform mixing time} of $\wt{\sfX}_N$ is defined as 
\begin{align}
\wt{\tau}(\eps) := \min \left\{ n \geq 0 : \left| \frac{ \wt{\sfP}^n(x,y) - \overline{\pi}_N(y) }{ \overline{\pi}_N(y) } \right| \leq \epsilon \text{ for all } x,y \in V_N^*  \right\}\,,
\label{eq:2_green_estimate_one_1}
\end{align}
where $\wt{\sfP}^n(x,y) := \wt{\sfP}_x( \wt{\sfX}_N(n) = y)$, and where $\overline{\pi}_N$ is the stationary measure associated to $\wt{\sfX}_N$, uniform on $V_N^*$. The discrete isoperimetric inequality on $V_N^*$ gives $a >0$, uniformly in $N$ for all $N$ large, such that for each nonempty subgraph $S \subset V_N^*$, we have $|\Delta S| \geq a | S|^{1/2}$. Here $\Delta S$ denotes the edge-boundary of $S$ in $V_N^*$. Using this bound in conjunction with Theorem 1 of \cite{morris2005evolving}, we obtain $\gamma \in \bbN$ so that $\wt{\tau}(1/2) \leq \gamma N^2$. 

Let $\wt{\sfG}_k^N$ be the analog of $\sfG_k^N$ for the lazy walk, and notice that 
\begin{align}
\wt{\sfG}_{\vart_N(n)}^N = 2 \sfG_{\vart_N(n)}^N \,.
\label{eq:lazy_central}
\end{align}
Decompose the trajectory of the lazy walk based on its location at time $\gamma N^2$:
\begin{align}
\wt{\sfG}_{\vart_N(n)}^N(0,x) = \wt{\sfG}^N_{\gamma N^2}(0,x) + \sum_{y \in V_N^*} \wt{\sfP}^{\gamma N^2} (0,y) \wt{\sfG}^N_{\vart_N(n) - \gamma N^2}(y, x) \,.
\label{eq:2_green_estimate_one_2}
\end{align}
By \eqref{eq:2_green_estimate_one_1}, the second term on the right side of \eqref{eq:2_green_estimate_one_2} obeys:
\begin{align}
\frac{1}{2} (n \log N - \gamma ) \leq \sum_{y \in V_N^*} \wt{\sfP}^{\gamma N^2} (0,y) \wt{\sfG}^N_{\vart_N(n) - \gamma N^2}(y, x) \leq \frac{3}{2} ( n \lceil \log N \rceil )\,,
\label{eq:2_green_estimate_one_3}
\end{align}
while by Lemma \ref{lem:step_green_estimate}, the first term on the right side of \eqref{eq:2_green_estimate_one_2} obeys
\begin{align}
0 \leq \sfG_{\gamma N^2}^N (0,x) \leq c \log N\,,
\label{eq:2_green_estimate_one_4}
\end{align}
where we have used Cauchy-Schwarz to bound $\sfG_{ \gamma N^2}^N (0,x)$ by $\sfG_{4 \gamma N^2}^N (0,0)$, for instance. Thus the constant $c$ in \eqref{eq:2_green_estimate_one_4} depends only on $\gamma$. We combine \eqref{eq:lazy_central}, \eqref{eq:2_green_estimate_one_3} and \eqref{eq:2_green_estimate_one_4}, taking $N$ larger if necessary to complete the proof. 
\end{proof}





\ACKNO{The work of A.C. and O.L. was supported in part by the European Union's - Seventh Framework Program (FP7/2007-2013) under grant agreement no. 276923 -- M-MOTIPROX. 
The work of A.C. was also supported by the Israeli Science Foundation grant no. 1723/14 -- Low Temperature Interfaces: Interactions, Fluctuations and Scaling; and by the Swiss National
Science Foundation 200021\underline{{ }{ }}163170. 
The work of O.L. was also supported by the Israeli Science Foundation grant no. 1328/17 and by grant I-2494-304.6/2017 from the German-Israeli Foundation
for Scientific Research and Development. 
The work of J.G. was supported in part by NSF grant DMS-1407558, by NSF grant DMS-1502632, and by a UCLA Dissertation Year Fellowship. }


\end{document}